\documentclass[a4paper,11pt]{amsart}
\usepackage{amssymb}
\usepackage{amsmath}
\usepackage{amsfonts}
\usepackage{color}

\newtheorem{X}{X}
\newtheorem{theorem}[X]{Theorem}
\theoremstyle{plain}

\newtheorem{lemma}[X]{Lemma}

\newtheorem{proposition}[X]{Proposition}

\numberwithin{equation}{section}
\input{tcilatex}

\begin{document}
\title[Dirichlet-Neumann problem]{On the explicit reconstruction of a
Riemann surface from its Dirichlet-Neumann operator}
\author{Gennadi Henkin}
\address[G. Henkin and V. Michel]{Universit\'{e} Pierre et Marie Curie\\
4, place Jussieu, 75252 Paris Cedex 05}
\email{henkin@math.jussieu.fr, michel@math.jussieu.fr}
\author{Vincent Michel}
\date{%
%TCIMACRO{\TeXButton{Date}{\today}}%
%BeginExpansion
\today%
%EndExpansion
}

\begin{abstract}
This article gives a complex analysis lighting on the problem which consists
in restoring a bordered connected riemaniann surface from its boundary and
its Dirichlet-Neumann operator. The three aspects of this problem, unicity,
reconstruction and characterization are approached.
\end{abstract}

\maketitle
\tableofcontents

\section{Statements of the main results}

Let $\mathcal{X}$ be an open bordered riemannian real surface (i.e. the
interior of an oriented riemannian two dimensional real manifold whose all
components have non trivial one dimensional smooth boundary) and $g$ its
metric. Using the boundary control method, Belishev and Kurylev~\cite%
{BeM1997}\cite{BeM-KuY1992} have started the study of the inverse problem
consisting in recovering $\left( \mathcal{X},g\right) $ from the operators $%
N_{\lambda }:C^{\infty }\left( b\mathcal{X}\right) \ni u\mapsto \left(
\partial \widetilde{u}_{\lambda }/\partial \nu \right) _{b\mathcal{X}}$
where $b\mathcal{X}$ is the boundary of $\mathcal{X}$, $\nu $ is the normal
exterior unit to $b\mathcal{X}$ and $\widetilde{u}_{\lambda }$ is the unique
solution of $\Delta _{g}U=\lambda U$ such that $U\left\vert _{b\mathcal{X}%
}\right. =u$. The principal result of \cite{BeM-KuY1992} implies that the
knowledge of $\lambda \mapsto N_{\lambda }$ on an non empty open set of $%
\mathbb{R}_{+}$ determines $\left( \mathcal{X},g\right) $ up to isometry.
The important question to know if $\left( \mathcal{X},g\right) $ is uniquely
determined by only one operator $N_{\lambda _{\ast }}$ with $\lambda _{\ast
}\neq 0$ remains open. This article mainly deals with the case of the
Dirichlet-Neumann operator $N_{\mathcal{X}}:=N_{0}$. Section~\ref{S/
Intrinsic} gives an intrinsic interpretation Electrical Impedance Tomography
on manifolds, EIT for short, in terms of Inverse-Dirichlet-Neumann problem
for twisted Laplacian. In dimension two, this clearly underline how the
complex structure of Riemannian surfaces is involved.\medskip

Two surfaces in the same conform class which have the same oriented boundary
and whose metrics coincides there, need to have the same Dirichlet-Neumann
operator. Conversely, Lassas and Uhlmann~\cite{LaM-UhG2001} have proved for
a connected $\mathcal{X}$ that the conform class and so the complex
structure of $\left( \mathcal{X},g\right) $ is determined by $N_{\mathcal{X}%
} $. Hence, it is relevant to consider $\mathcal{X}$ as a Riemann surface.
In~\cite{BeM2003}, using also the full knowledge of $N_{\mathcal{X}}$,
Belishev gives another proof of the above unicity by recovering abstractly $%
\mathcal{X}$ as the spectre of the algebra of boundary values of functions
holomorphic on $\mathcal{X}$ and continuous on $\overline{\mathcal{X}}=%
\mathcal{X}\cup b\mathcal{X}$. It turns out in our theorems~\ref{T/ unicite}
and ~\ref{T/ reconstr} that only three generic functions on the boundary and
their images by $N_{\mathcal{X}}$ are sufficient for unicity to hold and to
reconstruct $\mathcal{X}$ by integral Cauchy type formulas. Theorems~3a, 3b
and~3c deal with characterizations of data of the type $\left( b\mathcal{X}%
,N_{\mathcal{X}}\right) $ where $\mathcal{X}$ is a Riemann surface.\medskip

While the frame of bordered manifolds is sufficient for real analytic
boundaries, characterization statements lead to consider a wider class of
manifolds. In this article, $\left( \,\overline{\mathcal{X}},\gamma
\,\right) $ is a \textit{Riemann surface with almost smooth boundary} if the
following holds$^{\text{(}}$\footnote{%
A Stokes formula holds automatically for such manifolds (see lemma~\ref{L/
Stokes} in section~3).}$^{\text{)}}$~: $\overline{\mathcal{X}}$ is a compact
metrizable topological manifold which is the closure of $\mathcal{X}=%
\overline{\mathcal{X}}\backslash \gamma $, $\mathcal{X}$ is a Riemann surface%
$^{\text{(}}$\footnote{$h^{d}$ is the $d$-dimensionnial Hausdorff measure.}$%
^{\text{)}}$, $h^{2}\left( \,\overline{\mathcal{X}}\,\right) <\infty $, $%
\gamma $ is a smooth real curve$^{\text{(}}$\footnote{%
It could have been possible to allow singularities on $\gamma $ itself but
we have avoid it for the sake of simplicity of statements. Likewise, we
consider only smooth DN-datas in the sequel.}$^{\text{)}}$ and the set $%
\overline{\mathcal{X}}_{\func{sing}}$ of points of $\gamma $ where $%
\overline{\mathcal{X}}$ has not smooth boundary satisfies $h^{1}\left( 
\overline{\mathcal{X}}_{\func{sing}}\right) =0$~; $\overline{\mathcal{X}}%
\backslash \overline{\mathcal{X}}_{\func{sing}}$ is denoted $\overline{%
\mathcal{X}}_{\func{reg}}$.

If $\left( \,\overline{\mathcal{X}},\gamma \,\right) $ is a Riemann surface
with almost smooth boundary, classical results contained in \cite%
{AhL-SaL1960Livre} implies a Riemann's existence theorem~: a real valued
function $u$ of class $C^{1}$ on $\gamma $ has a unique continuous extension 
$\widetilde{u}$ to $\overline{\mathcal{X}}$ which is harmonic on $\mathcal{X}
$, smooth on $\overline{\mathcal{X}}_{\func{reg}}$ and satisfy $\int_{%
\mathcal{X}}i\,\partial \widetilde{u}\wedge \overline{\partial }\widetilde{u}%
<+\infty $. Moreover, $N_{\mathcal{X}}u$ still make sense as the element of
the dual space of $C^{1}\left( \gamma \right) $ which equals $\partial 
\widetilde{u}/\partial \nu $ on $\gamma \backslash \overline{\mathcal{X}}_{%
\func{sing}}$ (see proposition~\ref{L/ Di=Ri}).\medskip

In the sequel, $\gamma $ is a smooth compact oriented real curve without
component reduced to a point, $N$ is an operator from $C^{1}\left( \gamma
\right) $ to the space of currents on $\gamma $ of degree $0$ and order $1$
(i.e. functionals on $C^{1}$ 1-forms on $\gamma $), $\tau $ is a smooth
generating section of $T\gamma $ and $\nu $ is another vector field along $%
\gamma $ such that the bundle $\mathcal{T}$ generated by $\left( \nu
_{x},\tau _{x}\right) $, $x\in \gamma $, has rank $2$~; $\gamma $ is assumed
to be oriented by $\tau $ and $\mathcal{T}$ by $\left( \nu ,\tau \right) $.

The inverse Dirichlet-Neumann problem for $\left( \gamma ,N,\mathcal{T}%
\right) $ is to find, when it exists, an open riemaniann surface $\left( 
\mathcal{X},g\right) $ with almost smooth boundary $\gamma $ such that for
all $x\in \gamma \cap \overline{\mathcal{X}}_{\func{reg}}$, $\left( \nu
_{x},\tau _{x}\right) $ is a positively oriented orthonormal basis of $T_{x}%
\overline{\mathcal{X}}$ and for all $u\in C^{1}\left( \gamma \right) $, $%
Nu=N_{\mathcal{X}}u$ in the sense of currents. As these conditions do not
distinguished between metrics $g$ in a same conformal class, we look after $%
\mathcal{X}$ as a Riemann surface. The connection between real and complex
analysis in the IDN-problem is realized through the operators $L$ and $%
\theta $ defined for $u\in C^{1}\left( \gamma \right) $ by%
\begin{equation}
Lu=\frac{1}{2}\left( Nu-i\,Tu\right) \hspace{0.5cm}\&\hspace{0.5cm}\theta
u=\left( Lu\right) \left( \nu ^{\ast }+i\tau ^{\ast }\right)
\label{F/ N vers CR et theta}
\end{equation}%
where $T$ is the tangential derivation by $\tau $ and $(\nu _{x}^{\ast
},\tau _{x}^{\ast })$ is the dual basis of $\left( \nu _{x},\tau _{x}\right) 
$ for every $x\in \gamma $. Note that in the sense of currents, the equality 
$Nu=N_{\mathcal{X}}u$ is equivalent to the identity $\partial \widetilde{u}%
=\theta u$, the tilde denoting, as through all this article, continuous
harmonic extension to $\mathcal{X}$

If $\left( \overline{\mathcal{X}},\gamma \right) $ is a Riemann surface with
almost smooth boundary, $g$ is a hermitian metric on $\mathcal{X}$ for which 
$\left( \tau _{x},\nu _{x}\right) $ is a positively oriented orthonormal
basis of $T_{x}^{\ast }\mathcal{X}$ for $x\in \gamma $ outside $\sigma =%
\overline{\mathcal{X}}_{\func{sing}}$ and if $\rho \in C^{0}\left( \overline{%
\mathcal{X}}\right) \cap C^{\infty }\left( \overline{\mathcal{X}}\backslash
\sigma \right) $ is a defining function of $\gamma $ in $\overline{\mathcal{X%
}}$, then $(\nu ^{\ast },\tau ^{\ast })=\frac{1}{\left\vert d\rho
\right\vert _{g}}\left( d\rho ,d^{c}\rho \right) $ on $\gamma \backslash
\sigma $ where $d^{c}=i(\,\overline{\partial }-\partial \,)$ and $\partial 
\widetilde{u}=\left( Lu\right) \left\vert \partial \rho \right\vert
_{g}^{-1}\partial \rho =\theta u$ on $\gamma \backslash \sigma $ for all $%
u\in C^{1}\left( \gamma \right) $.\medskip

\noindent \textbf{Main hypothesis. }In addition to the assumptions on $%
\gamma $, we consider in all this paper $u_{0},u_{1},u_{2}\in C^{\infty
}\left( \gamma \right) $ three real valued functions only ruled by the main
hypothesis that%
\begin{equation}
f=\left( f_{1},f_{2}\right) =\left( \,\left( Lu_{\ell }\right) /\left(
Lu_{0}\right) \,\right) _{\ell =1,2}=\left( \,\left( \theta u_{\ell }\right)
/(\theta u_{0})\,\right) _{\ell =1,2}  \label{F/ f}
\end{equation}%
is an embedding of $\mathbb{\gamma }$ in $\mathbb{C}^{2}$ considered as the
complement of $\left\{ w_{0}=0\right\} $ in the complex projective plane $%
\mathbb{CP}_{2}$ with homogeneous coordinates $\left(
w_{0}:w_{1}:w_{2}\right) $. Prop.~\ref{L/ generique} whose proof is omitted
shows this is somehow generic~:

%TCIMACRO{\TeXButton{Compteur 0}{\setcounter{X}{-1}}}%
%BeginExpansion
\setcounter{X}{-1}%
%EndExpansion

\begin{proposition}
\label{L/ generique}Assume $\gamma $, $u_{0}$, $u_{1}$ real analytic and
that $f_{1}$ is non constant on each connected component of $\gamma $. For
any function $u_{2}\in C^{\omega }\left( \gamma \right) $, one can construct 
$v_{2}\in C^{\omega }\left( \gamma \right) $, arbitrarily close to $u_{2}$
in $C^{2}$ norm, such that $\left( f_{1},\left( Lv_{2}\right) /\left(
Lu_{0}\right) \,\right) $ is an embedding of $\mathbb{\gamma }$ into $%
\mathbb{C}^{2}$.
\end{proposition}

Assuming that $u=\left( u_{\ell }\right) _{0\leqslant \ell \leqslant 2}$
satisfies the main hypothesis, we set $\theta u=\left( \theta u_{\ell
}\right) _{0\leqslant \ell \leqslant 2}$ and call $\left( \gamma ,u,\theta
u\right) $ a \textit{restricted DN-datum} for an open Riemann surface $%
\mathcal{X}$ if $\mathcal{X}$ has almost smooth boundary $\gamma $, $\left(
\partial \widetilde{u_{\ell }}\right) {}\left\vert _{\gamma \backslash
\sigma }\right. =\theta u_{\ell }$ for $0\leqslant \ell \leqslant 2$ and the
well defined meromorphic quotient $F_{\ell }$ of (1,0)-forms $\left(
\partial \widetilde{u_{\ell }}\right) /\left( \partial \widetilde{u_{0}}%
\right) $ extends $f_{\ell }$ to $\mathcal{X}$ in the sense that for every $%
x_{0}\in \gamma $, $\underset{x\rightarrow x_{0},~x\in \mathcal{X}}{\lim }%
F\left( x\right) $ exists and equals $f\left( x_{0}\right) $. If $\gamma $
and $f$ are real analytic, this last property holds automatically.\medskip

We define an isomorphism between two Riemann surfaces with almost smooth
boundary $\left( \,\overline{\mathcal{X}},\gamma \,\right) $ and $\left( \,%
\overline{\mathcal{X}^{\prime }},\gamma ^{\prime }\,\right) $ as a map from $%
\overline{\mathcal{X}}$ to $\overline{\mathcal{X}^{\prime }}$ which realizes
a complex analytic isomorphism between $\mathcal{X}$ and $\mathcal{X}%
^{\prime }$. As the definition of a Riemann surface with almost smooth
boundary implies that its boundary is locally a Jordan curve in its double
which is the compact Riemann surface obtained by gluing along its boundary
its conjugate (see \cite{AhL-SaL1960Livre}), a theorem of Caratheodory
implies that if $\Phi :\mathcal{X}\rightarrow \mathcal{X}^{\prime }$ is a
complex analytic isomorphism, $\Phi $ and $\Phi ^{-1}$ extend continuously
to $\gamma $ and $\gamma ^{\prime }$ so that $\Phi $ becomes a homeomorphism
from $\overline{\mathcal{X}}$ to $\overline{\mathcal{X}^{\prime }}$. Hence, $%
\Phi $ is a diffeomorphism between manifolds with boundary from $\overline{%
\mathcal{X}}_{\func{reg}}\cap \Phi ^{-1}\left( \,\overline{\mathcal{X}%
^{\prime }}_{\func{reg}}\,\right) $ to $\overline{\mathcal{X}^{\prime }}_{%
\func{reg}}\cap \Phi \left( \overline{\mathcal{X}}_{\func{reg}}\right) $.

The first theorem of this article is a significative improvement of results
in \cite{BeM-KuY1992}\cite{LaM-UhG2001} on how unique $\mathcal{X}$ can be
when a restricted DN-datum is specified.

\begin{theorem}
\label{T/ unicite}Assume that $\mathcal{X}$ and $\mathcal{X}^{\prime }$ are
open Riemann surfaces with \textit{restricted DN-datum} $\left( \gamma
,u,\theta u\right) $. Then, there is an isomorphism of Riemann surfaces with
almost smooth boundary between $\mathcal{X}\cup \gamma $ and $\mathcal{X}%
^{\prime }\cup \gamma $ whose restriction on $\gamma $ is the identity.
\end{theorem}

\noindent \textbf{Remarks. 1. }If $E\subset \gamma $ and $h^{1}\left( E\cap
c\right) >0$ for each connected component $c$ of $\gamma $, meromorphic
functions are uniquely determined by their values on $E$ and it follows that
th.~\ref{T/ unicite} conclusions hold when $N_{\mathcal{X}^{\prime }}u_{\ell
}=N_{\mathcal{X}}u_{\ell }$ is ensured only on $E$ and the meromorphic
functions $\left( \partial \widetilde{u_{\ell }}\right) /\left( \partial 
\widetilde{u}_{0}\right) $ are continuous near $\gamma $. This includes \cite%
[th. 1.1.i]{LaM-UhG2001} which is stated for a connected $\mathcal{X}$%
.\smallskip

\textbf{2. }The proof of theorem~\ref{T/ unicite} also contains the fact
that two connected compact Riemann surfaces $\mathcal{Z}$ and $\mathcal{Z}%
^{\prime }$ are isomorphic when they share a same real smooth curve $\gamma $
which can be embedded into $\mathbb{C}^{2}$ by a map which extends
meromorphically both to $\mathcal{Z}$ and $\mathcal{Z}^{\prime }$ and
continuously near $\gamma $.\medskip

The assumption on $u_{0}$, $u_{1}$ and $u_{2}$ is used only to ensure that
the map $f$ defined by (\ref{F/ f}) is an embedding of $\gamma $ into $%
\mathbb{C}^{2}$ extending meromorphically to $\mathcal{X}$ into $F=\left(
\partial \widetilde{u_{\ell }}/\partial \widetilde{u}_{0}\right) _{\ell
=1,2} $. Moreover, theorem~\ref{T/ unicite de Y} in section~\ref{S/ UsousE}
implies that if $\mathcal{X}$ has almost smooth boundary and solves the
IDN-problem, the map $F$ enable to see $\overline{\mathcal{X}}$ as a
normalization of the closure of a complex curve$^{\text{(}}$\footnote{%
That is a pure 1-dimensionnal complex analytic subset of an open subset of $%
\mathbb{CP}_{2}$.}$^{\text{)}}$ of $\mathbb{CP}_{2}\backslash f\left( \gamma
\right) $ uniquely determined by $\gamma $. This shows that in each
characterization theorem~\ref{T/ caract - a}, \ref{T/ caract - b}, \ref{T/
caract - c}, the constructed Riemann surface is, up to isomorphism, the only
one which has a chance to solve the IDN-problem.\medskip

Our next result explains how to recover $F\left( \mathcal{X}\right) $ and $%
\partial \widetilde{u_{\ell }}$ from $\theta u_{\ell }$ and the intersection
of $F\left( \mathcal{X}\right) $ with the lines $\Delta _{\xi }=\{z_{2}:=%
\frac{w_{2}}{w_{0}}=\xi \}$, $\xi \in \mathbb{C}$. Desingularization
arguments enable then the reconstruction of $\mathcal{X}$ from $F\left( 
\mathcal{X}\right) $.

\begin{theorem}
\label{T/ reconstr}If $\mathcal{X}$ is an open Riemann surface with \textit{%
restricted DN-datum} $\left( \gamma ,u,\theta u\right) $, the following
holds~:

1) The map $f$ defined by (\ref{F/ f}) has a meromorphic extension $F$ to $%
\mathcal{X}$ and there are discrete sets $\mathcal{A}$ and $\mathcal{B}$ in $%
\mathcal{X}$ and $\mathcal{Y}=F\left( \mathcal{X}\right) \backslash f\left(
\gamma \right) $ respectively such that $F:\mathcal{X}\backslash \mathcal{A}%
\rightarrow \mathcal{Y}\backslash \mathcal{B}$ is one to one.

2) Almost all $\xi _{\ast }\in \mathbb{C}$ has a neighborhood $W_{\xi _{\ast
}}$ such that for all $\xi $ in $W_{\xi _{\ast }}$, $\mathcal{Y}_{\xi }=%
\mathcal{Y}\cap \Delta _{\xi }=\underset{1\leqslant j\leqslant p}{\cup }%
\{(h_{j}\left( \xi \right) ,\xi )\}$ where $h_{1},...,h_{p}$ are $p$
mutually distinct holomorphic functions on $W_{\xi _{\ast }}$ whose
symmetric functions $S_{h,m}=\underset{1\leqslant j\leqslant p}{\Sigma }%
h_{j}^{m}$ are recovered by the Cauchy type integral formulas%
\begin{equation}
\frac{1}{2\pi i}\int_{\gamma }\frac{f_{1}^{m}}{f_{2}-\xi }%
df_{2}=S_{h,m}\left( \xi \right) +P_{m}\left( \xi \right) ,~m\in \mathbb{N},
\tag{\TeXButton{E_m_ksi}{$E_{m,\protect\xi
}$}}  \label{F/ equation pour Y}
\end{equation}%
where $P_{m}$ is a polynomial of degree at most $m$. More precisely, the
system $E_{\xi }=(E_{m,\xi _{\nu }})_{\substack{ 0\leqslant m\leqslant B-1 
\\ o\leqslant \nu \leqslant A-1}}$ enables explicit computation of $%
h_{j}\left( \xi _{\nu }\right) $ and $P_{m}$ if $A\geqslant B\geqslant 2p+1$
and $\xi _{0},...,\xi _{A}$ are mutually distinct points.

3) For almost all $\xi _{\ast }\in \mathbb{C}$, $W_{\xi _{\ast }}$ can be
chosen so that $\mathcal{B}\cap \underset{\xi \in W_{\xi _{\ast }}}{\cup }%
\mathcal{Y}_{\xi }=\varnothing $ and $\partial \widetilde{u_{\ell }}$, $%
0\leqslant \ell \leqslant 2$, can be reconstructed in $F^{-1}(\underset{\xi
\in W_{\xi _{\ast }}}{\cup }\mathcal{Y}_{\xi })$ from the well defined
meromorphic quotient $\left( \partial \widetilde{u_{\ell }}\right) /\left(
\partial F_{2}\right) $ thanks to the Cauchy type formulas 
\begin{equation}
\frac{1}{2\pi i}\int_{\gamma }\frac{f_{1}^{m}}{f_{2}-\xi }\theta u_{\ell
}=\dsum\limits_{1\leqslant j\leqslant p}h_{j}\left( \xi \right) ^{m}\frac{%
\partial \widetilde{u_{\ell }}}{\partial F_{2}}\left( F^{-1}\left(
h_{j}\left( \xi \right) ,\xi \right) \right) \mathbb{+}Q_{m}\left( \xi
\right)  \tag{\TeXButton{E_tilde_m_ksi}{$T_{m,\protect\xi
}$}}
\label{F/ equation Theta}
\end{equation}%
where $m$ is any integer and $Q$ is a polynomial of degree at most $m$.
\end{theorem}

\noindent \textbf{Remark. }The number $\alpha $ of connected components of $%
\mathcal{X}$ can be computed by the following algorithm~: let $\gamma _{1}$
be a component of $\gamma $ and let $\lambda _{1}$ be a function which is
zero on $\gamma _{j}$ for $j\neq 0$ and non constant on $\gamma _{1}$~; then
if $\mathcal{X}_{1}$ is the component of $\mathcal{X}$ whose boundary
contains $\gamma _{1}$, $N\lambda _{1}\neq 0$ on each component $\gamma
_{1},...,\gamma _{k}$ of $\gamma $ which with $\gamma _{1}$ are the
components of $b\mathcal{X}_{1}$. Iterating this with components $\gamma $
different from $\gamma _{1},...,\gamma _{k}$, yields a process with $\alpha $
steps.\smallskip

The numerical resolution of $\left( E_{\xi }\right) $ and the study of its
stability requires an estimate of the number $I_{\Delta _{\xi }}$ of points
of intersection, multiplicities taken in account, of $\mathcal{Y}$ with $%
\Delta _{\xi }$. To achieve this, it is sufficient to estimates the number $%
I_{\Delta }$ of intersection points of $\mathcal{Y}$ with a $\mathbb{CP}_{2}$%
-line $\Delta $ generic in the sense that $\Delta $ does not contained the
germ of a component of $\mathcal{Y}$ near $\gamma $. Indeed, if $L$ (resp. $%
L_{\xi }$) denotes a linear homogeneous form defining $\Delta $ (resp. $%
\Delta _{\xi }$),\vspace{-2pt}%
\begin{equation*}
I_{\Delta _{\xi }}-I_{\Delta }=\frac{1}{2\pi i}\int_{\gamma }\left( L_{\xi
}/L\right) ^{-1}d\left( L_{\xi }/L\right) .\vspace{-2pt}
\end{equation*}%
Thus, an a priori upper bound of $I_{\Delta }$ for any particular line $%
\Delta $ would be very useful. This open problem is related, because of the
Ahlfors theorem on covering surface, to the computation of the genus $g_{%
\mathcal{X}}$ of $\mathcal{X}$ from some DN-datum when $\mathcal{X}$ is
connected. Under the condition $\gamma $ is connected, Belishev~\cite%
{BeM2003} has shown that $2g_{\mathcal{X}}$ is the rank of $Id+\left( N_{%
\mathcal{X}}T^{-1}\right) ^{2}$ acting on the space of smooth functions on $%
\gamma $ admitting a smooth primitive, $T^{-1}$ being a primitive operator.
Yet, a formula for $g_{\mathcal{X}}$ involving only the action of $N_{%
\mathcal{X}}$ on a finite generic set of functions has to be found.\smallskip

The third aspect of the IDN-problem, characterization of what can and should
be a DN-datum, have lead us to allow $\mathcal{X}$ to have only almost
smooth boundary. Th.~\ref{T/ caract - a} below explicitly characterizes the
only right candidate for $\mathcal{X}$ while its part~C gives a test which
discriminates which $\left( \gamma ,u,\theta u\right) $ are DN-data and
which are not. To perform it, we need a \textit{Green function for} $%
\overline{\mathcal{X}}$ relatively to a domain $\mathcal{D}$ of $\mathcal{Z}$
containing $\overline{\mathcal{X}}$, that is a smooth symmetric function $g$
defined on $\mathcal{D}\times \mathcal{D}$ without its diagonal such that
each $g\left( .,z\right) $ is harmonic on $\mathcal{D}\backslash \left\{
z\right\} $ and has singularity $\frac{1}{2\pi }\ln \func{dist}\left(
.,z\right) $ at $z$, the distance being computed in any hermitian metric on $%
\mathcal{Z}$.

%TCIMACRO{%
%\TeXButton{Compteurs spéciaux}{\newcounter{Y}\setcounter{Y}{1}
%\renewcommand{\theX}{\arabic{X}\alph{Y}}}}%
%BeginExpansion
\newcounter{Y}\setcounter{Y}{1}
\renewcommand{\theX}{\arabic{X}\alph{Y}}%
%EndExpansion

\begin{theorem}
\label{T/ caract - a}Assume that the main hypothesis is valid and consider$%
\vspace{-3pt}$%
\begin{equation}
G:\mathbb{C}^{2}\ni \left( \xi _{0},\xi _{1}\right) \mapsto \frac{1}{2\pi i}%
\int_{\gamma }f_{1}\frac{d\left( \xi _{0}+\xi _{1}f_{1}+f_{2}\right) }{\xi
_{0}+\xi _{1}f_{1}+f_{2}}.\vspace{-3pt}  \label{F/ fonction G}
\end{equation}

\textbf{A.} If an open Riemann surface $\mathcal{X}$ has restricted DN-datum 
$\left( \gamma ,u,\theta u\right) $, then almost all point $\xi _{\ast }$\
of $\mathbb{C}^{2}$\ has a neighborhood where one can find mutually distinct
holomorphic functions $h_{1},...,h_{p}$\ such that%
\begin{gather}
0=\frac{\partial ^{2}}{\partial \xi _{0}^{2}}(G-\dsum\limits_{1\leqslant
j\leqslant p}h_{j})  \label{F/ caract G} \\
h_{j}\frac{\partial h_{j}}{\partial \xi _{0}}=\frac{\partial h_{j}}{\partial
\xi _{1}},1\leqslant j\leqslant p.  \label{F/ SW}
\end{gather}

\textbf{B.} Conversely, assume $\gamma $ is connected and the conclusion of
A is satisfied in a connected neighborhood $W_{\xi _{\ast }}$ of one point $%
\left( \xi _{0\ast },\xi _{1\ast }\right) $. Then, if $\left( \partial
^{2}G/\partial \xi _{0}^{2}\right) _{\left\vert W_{\xi _{\ast }}\right.
}\neq 0$, there is an open Riemann surface $\mathcal{X}$ with almost smooth
boundary$^{\text{(}}$\footnote{%
With \cite[example 10.5]{HaR-LaB1975}, one can construct smooth restricted
DN-datas for which the solution of the IDN-problem is a manifold with only
almost smooth boundary.}$^{\text{)}}$ $\gamma $ where $f$\ extends
meromorphically. If $\left( \partial ^{2}G/\partial \xi _{0}^{2}\right)
_{\left\vert W_{\xi _{\ast }}\right. }=0$, the same conclusion holds for a
suitable orientation of $\gamma $.

\textbf{C.} Assume that $\left( \,\overline{\mathcal{X}},\gamma \,\right) $\
is a Riemann surface with almost smooth boundary. Let $\mathcal{Z}$\ be the
double of $\mathcal{X}$, $\mathcal{D}$ a smooth domain of $\mathcal{Z}$
containing $\overline{\mathcal{X}}$ and $g$ a Green function for $\overline{%
\mathcal{X}}$ relatively to $\mathcal{D}$. Then, $\left( \gamma ,u,\theta
u\right) $\ is actually a restricted DN-datum if and only if for any $z\in 
\mathcal{D}\backslash \overline{\mathcal{X}}$,%
\begin{equation}
\int_{\gamma }u_{\ell }\left( \zeta \right) \partial _{\zeta }g\left( \zeta
,z\right) +g\left( \zeta ,z\right) \overline{\theta u_{\ell }\left( \zeta
\right) }=0,~0\leqslant \ell \leqslant 2.  \label{F/ Green}
\end{equation}
\end{theorem}

\noindent \textbf{Remarks. 1. }The connectness of $\gamma $ is essentially
used to ensure that any possible solution to the IDN-problem has to be
connected. Taking in account the remark following th.~\ref{T/ reconstr}, one
may weaken the connectness assumption on $\gamma $ into the requirement that
the given DN-datum ensures that possible solutions are connected. Then, the
conclusions of th.~\ref{T/ caract - a}.B are still true (see the proof).

\textbf{2. }The proof includes that if $\gamma $ and $f$ are real analytic, $%
\left( \,\overline{\mathcal{X}},\gamma \,\right) $ is a manifold with
boundary in the classical sense.

\textbf{3. }Emphasizing on $f_{2}$ instead of $f_{1}$, one can consider%
\begin{equation*}
G_{2}:\xi \mapsto \frac{1}{2\pi i}\int_{\gamma }f_{2}\frac{d\left( \xi
_{0}+\xi _{1}f_{1}+f_{2}\right) }{\xi _{0}+\xi _{1}f_{1}+f_{2}}.
\end{equation*}%
If $h_{j}$ is linked to $h_{j,2}$ by $0=\xi _{0}+\xi _{1}h_{j}+h_{j,2}$, $%
\left( h_{1},...,h_{p}\right) $ satisfy (\ref{F/ SW}) and (\ref{F/ caract G}%
) if and only $\frac{\partial ^{2}}{\partial \xi _{0}^{2}}(G_{2}-\underset{%
1\leqslant j\leqslant p}{\Sigma }h_{j,2})=0$ and $h_{j}\frac{\partial h_{j,2}%
}{\partial \xi _{0}}=\frac{\partial h_{j,2}}{\partial \xi _{1}}$, $%
1\leqslant j\leqslant p$.

\textbf{4. }Select $\mathcal{H}=\left\{ h_{1},...,h_{p}\right\} $ satisfying
(\ref{F/ SW}) and minimal for (\ref{F/ caract G}). Then, section~\ref{S/ 3aB}
and proposition~\ref{L/ caract sol} shows that there is $\tau ^{\prime
}\subset \delta =f\left( \gamma \right) $ such that $h^{1}\left( \tau
^{\prime }\right) =0$ and $\mathcal{X}$ is a normalization of the abstract
curve $\mathcal{Y}\cup \tau ^{\prime }$ where, when $\mathcal{H}=\varnothing 
$, $\mathcal{Y}$ is the polynomial hull of $\delta $ in the affine complex
plane 
\begin{equation*}
\mathbb{C}_{\xi _{\ast }}^{2}=\left\{ w\in \mathbb{CP}^{2}~;~\xi _{\ast
}w=\xi _{0\ast }w_{0}+\xi _{1\ast }w_{1}+w_{2}\neq 0\right\} .
\end{equation*}%
and, otherwise, $\mathcal{Y}$ is the analytic extension in $\mathbb{CP}%
_{2}\backslash \delta $ of the union of the graphs of the functions $\left(
1:h_{j}:-\xi _{0}-\xi _{1}h_{j}\right) $, $1\leqslant j\leqslant p$. Hence,
when $\mathcal{H}$ is minimal, decomposition (\ref{F/ caract G}) of $G$ is
unique up to order and $\func{Card}\mathcal{H}$ is the minimal number $p$
for which such a decomposition exists. Moreover, theorem~\ref{T/ unicite de
Y} implies that the only Riemann surfaces $\mathcal{X}$ which has a chance
to solve the IDN-problem are normalizations of $\mathcal{Y}$.$\medskip $

The vanishing of $\partial ^{2}G/\partial \xi _{0}^{2}$ in a connected
neighborhood $W_{\xi _{\ast }}$ of $\xi _{\ast }\in \mathbb{C}^{2}$ is known
to be equivalent to the fact that $\delta =f\left( \gamma \right) $ satisfy
the classical Wermer-Harvey-Lawson moment condition in$^{\text{(}}$\footnote{%
When $\xi _{\ast }$ belongs to the connected component of infinity of $%
\left\{ \xi \in \mathbb{C}^{2}~;~\forall w\in \delta ,~\xi w\neq 0\right\} $%
, this moment condition is equivalent to the moment condition in $\mathbb{C}%
_{(1,0)}^{2}$ and also to the vanishing of $G$ on this component.}$^{\text{)}%
}$ $\mathbb{C}_{\xi _{\ast }}^{2}~$: for all $k_{1},k_{2}\in \mathbb{N}$, $%
\int_{\delta }z_{1}^{k_{1}}z_{2}^{k_{2}}dz_{2}=0$ where $z=\left( w_{j}/\xi
_{\ast }w\right) _{j=1,2}$ (see \cite[cor. 1.6.2]{DoP-HeG1997}). It is
proved in \cite{WeJ1958c} for the real analytic case and in \cite{BiE1963}%
\cite{HaR-LaB1975} for the smooth case that for a suitable orientation of $%
\gamma $, this moment condition guarantees the existence in $\mathbb{C}_{\xi
_{\ast }}^{2}\backslash \delta $ of a unique complex curve $\mathcal{Y}$
with finite mass and boundary $\pm \delta $ in the sense of currents. In 
\cite{AlH-WeJ2000}, Alexander and Wermer have improved this
Wermer-Bishop-Harvey-Lawson statement by showing that a closed oriented
smooth connected real curve $\delta $ of $\mathbb{C}^{2}$ is, with its given
orientation, the boundary, in the sense of currents, of a complex curve of
finite mass in $\mathbb{C}_{\xi _{\ast }}^{2}\backslash \delta $ if and only
if $\frac{1}{2\pi i}\int_{\delta }\frac{dA}{A}\geqslant 0$ for any
polynomial $A$ which does not vanish on $\gamma $. Hence, in case $\left(
\partial ^{2}G/\partial \xi _{0}^{2}\right) _{\left\vert W_{\xi _{\ast
}}\right. }=0$ it is sufficient to find one polynomial $A$ such that $%
\int_{f\left( \gamma \right) }\frac{dA}{A}\neq 0$ to determine the correct
orientation of $\gamma $.

Note that the case $\left( \partial ^{2}G/\partial \xi _{0}^{2}\right)
_{\left\vert W_{\xi _{\ast }}\right. }=0$ occurs only for very special
DN-data since it implies that for $\ell =1,2$, $f_{\ell }$ admits a $\mathbb{%
C}_{\xi _{\ast }}^{2}$-valued holomorphic extension to $\mathcal{X}$.
Proposition~\ref{T/ affine} proposes another result of this kind for some
other special DN-data when they are available.\medskip

To palliate the difficulty of computing Green functions, the theorem below
proposes another way to achieve the same goals~: select the right candidates
for $\mathcal{X}$ and extension of $\theta u_{\ell }$; check this yields a
solution.

%TCIMACRO{\TeXButton{+1/-1}{\addtocounter{X}{-1}\stepcounter{Y}}}%
%BeginExpansion
\addtocounter{X}{-1}\stepcounter{Y}%
%EndExpansion

\begin{theorem}
\label{T/ caract - b}Assume that the main hypothesis is valid. Let $G$ be
the function defined by (\ref{F/ fonction G}) and let be $\widetilde{G}$ the
form which in $\mathbb{CP}_{2}$ with homogenous coordinates $\eta =(\eta
_{0}:\eta _{1}:\eta _{2})$ is given by%
\begin{equation}
\widetilde{G}=\dsum\limits_{0\leqslant \ell \leqslant 2}\frac{1}{2\pi i}%
\left( \int_{\gamma }\frac{\theta u_{\ell }}{\eta _{0}+\eta _{1}f_{1}+\eta
_{2}f_{2}}\right) d\eta _{\ell }=\dsum\limits_{0\leqslant \ell \leqslant 2}%
\widetilde{G}_{\ell }\,d\eta _{\ell }  \label{F/ formes G}
\end{equation}

\textbf{A. }If an open Riemann surface $\mathcal{X}$ has restricted DN-datum 
$\left( \gamma ,u,\theta u\right) $, then,

\noindent \textbf{a1) }Almost all points $\eta _{\ast }=\left( \xi _{\ast
0}:\xi _{\ast 1}:1\right) $ of $\mathbb{CP}_{2}$ has a neighborhood where $%
\widetilde{G}$ can be written as the sum of $p$ holomorphic closed forms $%
g_{j}=\underset{0\leqslant \ell \leqslant 2}{\Sigma }g_{j,\ell }\,d\eta
_{\ell }$ such that $\left( h_{j}\right) =\left( g_{j,1}/g_{j,0}\right)
_{1\leqslant j\leqslant p}$ satisfy (\ref{F/ SW}).

\noindent \textbf{a2) }The form $\Theta _{\ell }=\partial \widetilde{u_{\ell
}}$, $0\leqslant \ell \leqslant 2$, satisfies 
\begin{equation}
\int_{c}\func{Re}\Theta _{\ell }=0  \label{F/ sans Green1}
\end{equation}%
for all $c$ in the first homology group $H_{1}\left( \mathcal{X}\right) $ of 
$\mathcal{X}$.

\textbf{B. b1) }Assume $\gamma $ is connected and there is $\eta _{\ast
}=\left( \xi _{\ast 0}:\xi _{\ast 1}:1\right) $ and a connected neighborhood 
$W_{\xi _{\ast }}$ of $\xi _{\ast }=\left( \xi _{\ast 0},\xi _{\ast
1}\right) $ such that (a1) is true for all $\eta \in W_{\eta _{\ast
}}=\left\{ \left( \xi _{0}:\xi _{1}:1\right) ~;~\left( \xi _{0},\xi
_{1}\right) \in W_{\xi _{\ast }}\right\} $. Then, there exists an open
Riemann surface $\mathcal{X}$, topologically bordered by $\gamma $, where $f$
extends meromorphically and each $\theta u_{\ell }$ extends weakly$^{\text{(}%
}$\footnote{%
It means that $\int_{\gamma }\varphi \theta u_{\ell }=\int_{\mathcal{X}%
}d\left( \varphi \Theta _{\ell }\right) =\int_{\mathcal{X}}\left( \,%
\overline{\partial }\varphi \,\right) \wedge \Theta _{\ell }$ holds for any
Lipschitz function $\varphi $ on $\overline{\mathcal{X}}$ which is a
holomorphic function of $f$ near points of $\Sigma $ and singular points of $%
\left( \,\overline{\mathcal{X}},\gamma \,\right) $~; if $\left( \,\overline{%
\mathcal{X}},\gamma \,\right) $ is a manifold with boundary, this definition
means that $\Theta _{\ell }\left\vert _{\gamma }\right. =\theta u_{\ell }$
in the usual sense.}$^{\text{)}}$ into a meromorphic (1,0)-form $\Theta
_{\ell }$ outside a set $\Sigma $ of zero length$^{\text{(}}$\footnote{%
Basing on \cite[example 1]{DiT1998b}, one can construct examples where (a1)
is satisfied while the weak extension $\Theta _{\ell }$ has essential
singularities on some zero length set $\Sigma $.}$^{\text{)}}$.

\noindent \textbf{b2) }In addition to (a1), assume that $\Theta _{\ell }$
satisfies $\int i\,\Theta _{\ell }\wedge \overline{\Theta _{\ell }}<+\infty $
and (\ref{F/ sans Green1}). Then if $\left( \partial ^{2}G/\partial \xi
_{0}^{2}\right) _{\left\vert W_{\xi _{\ast }}\right. }\neq 0$, $\left( \,%
\overline{\mathcal{X}},\gamma \,\right) $ is a manifold with almost smooth
boundary~; the same conclusion holds when $\widetilde{G}\left\vert _{W_{\eta
_{\ast }}}\right. =0$ if $\gamma $ has a suitable orientation. If $\left(
\partial ^{2}G/\partial \xi _{0}^{2}\right) _{\left\vert W_{\xi _{\ast
}}\right. }=0$ but $\widetilde{G}\left\vert _{W_{\eta _{\ast }}}\right. \neq
0$, then either $\mathcal{X}$ is a domain with boundary $\gamma $ in a
normalization of an algebraic curve of $\mathbb{CP}_{2}$, either $\overline{%
\mathcal{X}}$ is a compact Riemann surface where $\gamma $ is a slit$^{\text{%
(}}$\footnote{%
This means that $\overline{\mathcal{X}}\backslash \gamma $ is connected.}$^{%
\text{)}}$. In all cases, $u_{\ell }$ admits a continuous extension $%
\widetilde{u_{\ell }}$ to $\overline{\mathcal{X}}$ which is harmonic in $%
\mathcal{X}$ and such that $\Theta _{\ell }=\partial \widetilde{u_{\ell }}$,
which means that $Nu_{\ell }$ is actually the DN-datum of $\mathcal{X}$ for $%
u_{\ell }$.
\end{theorem}

\noindent \textbf{Remarks. 1}. When (a1) holds, $h_{j,2}=\frac{g_{j,2}}{%
g_{j,0}}$ verify $h_{j}\frac{\partial h_{j,2}}{\partial \xi _{0}}=\frac{%
\partial h_{j,2}}{\partial \xi _{1}}$, $1\leqslant j\leqslant p$.

\textbf{2. }Formulas~(\ref{F/ equation pour Y}) and~(\ref{F/ equation Theta}%
) enable direct reconstruction of a projective presentation of $\mathcal{X}$
and forms $\Theta _{\ell }$.\medskip

Theorem~\ref{T/ caract - b} is obtained by a normalization of a singular
version of the IDN-problem which is more explicit. When $\mathcal{X}$ is
smooth, the harmonicity of a distribution $U$ is equivalent to the fact that 
$\partial U$ is holomorphic. For the case where $\mathcal{X}$ is a complex
curve of an open set in $\mathbb{CP}_{2}$, we need two of the several non
equivalent definitions of holomorphic (1,0)-forms.

At first, we use the \textit{weakly holomorphic} forms introduced by
Rosenlicht~\cite{RoM1952} which can be defined as meromorphic (1,0)-forms $%
\psi $ such that $\psi \wedge \lbrack \mathcal{X]}$ is $\overline{\partial }$%
-closed current of $\mathbb{CP}_{2}$. Such forms $\psi $ are also
characterized by the fact that $p_{\ast }\psi $ is a usual holomorphic $%
\left( 1,0\right) $-form for any holomorphic proper function $p:\mathcal{X}%
\rightarrow \mathbb{C}$. A distribution $U$ is defined as \textit{weakly
harmonic} if $\partial U$ is weakly holomorphic.

Assume now $\mathcal{X}$ lies in $\mathbb{CP}_{2}$ and that $\mathcal{X}$ is
bounded in the sense of currents by $\gamma $. A distribution $U$ on $%
\mathcal{X}$ is said almost smooth up to the boundary if it is the case near
each $p\in \gamma $ where $\left( \,\overline{\mathcal{X}},\gamma \,\right) $
is a manifold with boundary and if $U$ has a restriction on $\gamma $ in the
sense of currents. When $u$ is a smooth function on $\gamma $, a weakly
harmonic extension of $u$ to $\mathcal{X}$ is a weakly harmonic distribution 
$U$ almost smooth up to boundary whose restriction on $\gamma $ is $u$.
Since two weakly harmonic extension $U_{1}$ and $U_{2}$ of $u$ to $\mathcal{X%
}$ are equal when $\partial U_{1}=\partial U_{2}$ on $\gamma $ in the sense
of currents, we consider a \textit{weak} \textit{Cauchy-Dirichlet problem}~:
a data is a smooth function $u$ on $\gamma $ and a smooth section $\lambda $
of $T_{\gamma }^{\ast }\mathcal{X~}$; a solution is a weakly harmonic
function $U$ almost smooth up to $\gamma $ such that $u=U\left\vert _{\gamma
}\right. $ and $\lambda =\left( \partial U\right) _{\gamma }$ in the sense
of currents~; when it exists, such an $U$ is unique and is denoted $%
\widetilde{u}$ as any harmonic extension in this article. In connection with
this notion, we define a \textit{weak restricted data }as a triplet $%
(\,\gamma ,u,\theta u\,)$ where $u=\left( u_{\ell }\right) _{0\leqslant \ell
\leqslant 2}$ (resp. $\theta u=\left( \theta u_{\ell }\right) _{0\leqslant
\ell \leqslant 2}$) is a triplet of smooth functions (resp. (1,0)-forms) on $%
\gamma $ such that $\theta u_{\ell }=\left( \partial \widetilde{u_{\ell }}%
\right) _{\gamma }$ in the sense of currents.\smallskip

The weak CD-problem has its own interest and arise naturally in the proof of
th.~\ref{T/ caract - b}. However, the original IDN-problem requires a more
restrictive notion of harmonicity. According to Griffiths~\cite{GrP1976}, 
\textit{holomorphic forms} (resp. \textit{harmonic functions}) are by
definition push forward of holomorphic forms (resp. harmonic functions) on a
normalization of $\mathcal{X}$. Equivalently, a real function $U$ on $%
\mathcal{X}$ is harmonic if and only if $U$ is harmonic in the regular part $%
\mathcal{X}_{\func{reg}}$ of $\mathcal{X}$ and $\int_{\mathcal{X}_{\func{reg}%
}}i\,\partial U\wedge \overline{\partial }U<+\infty $. This notion is close
in spirit to a Riemann characterization of the harmonic function with given
boundary value $u$ as the smooth function extending $u$ to $\overline{%
\mathcal{X}}$ and minimizing the preceding integral. We can now state a
singular version of th.~\ref{T/ caract - b}.

%TCIMACRO{\TeXButton{+1/-1}{\addtocounter{X}{-1}\stepcounter{Y}}}%
%BeginExpansion
\addtocounter{X}{-1}\stepcounter{Y}%
%EndExpansion

\begin{theorem}
\label{T/ caract - c}Consider in $\mathbb{CP}_{2}\backslash \left\{
w_{0}=0\right\} $ a smooth oriented real curve $\gamma $, three functions $%
u_{0}$, $u_{1}$, $u_{2}$ in $C^{\infty }\left( \gamma \right) $ and $\theta
_{0}$, $\theta _{1}$, $\theta _{2}$ three smooth sections of $\Lambda
^{1,0}T_{\gamma }^{\ast }\mathbb{CP}_{2}$ such that $du_{\ell }=2\func{Re}%
\theta _{\ell }$, $0\leqslant \ell \leqslant 2$ and linked by the relations $%
\theta _{1}=z_{1}\theta _{0}$, $\theta _{2}=z_{2}\theta _{0}$. Let $G$ and $%
\widetilde{G}$ be the form given by (\ref{F/ fonction G}) and (\ref{F/
formes G}) but with $\left( f_{1},f_{2}\right) =(z_{1},z_{2})$.

\textbf{A. }Assume $\gamma $ bounds, in the sense of currents, a complex
curve $\mathcal{X}$ of $\mathbb{CP}_{2}\backslash \gamma $ which has finite
volume and weak restricted DN-datum $\left( \gamma ,u,\theta u\right) $.
Then,

\noindent \textbf{a1) }The conclusions of theorem~\ref{T/ caract - b}.A.a1
are valid.

\noindent \textbf{a2) }The form $\Theta _{\ell }=\partial \widetilde{u_{\ell
}}$ satisfies (\ref{F/ sans Green1}) for all $c$ in $H_{1}\left( \mathcal{X}%
_{\func{reg}}\right) $.

\textbf{B. b1) }Conversely, assume that $\gamma $ is connected and that (a1)
is valid for one point $\eta _{\ast }=\left( \xi _{0\ast }:\xi _{1\ast
}:1\right) $. Then, there is in $\mathbb{CP}_{2}\backslash \gamma $ a
complex curve $\mathcal{X}$ of finite mass where each $\theta _{\ell }$
extends weakly on $\mathcal{X}$ into a weakly holomorphic $\left( 1,0\right) 
$-form $\Theta _{\ell }$. Moreover, if $\left( \partial ^{2}G/\partial \xi
_{0}^{2}\right) _{\left\vert W_{\xi _{\ast }}\right. }\neq 0$, then $%
\mathcal{X}$ has boundary $\gamma $ in the sense of currents~; the same
conclusion holds if $\widetilde{G}\left\vert _{W_{\eta _{\ast }}}\right. =0$
but for a suitable orientation of $\gamma $. If $\left( \partial
^{2}G/\partial \xi _{0}^{2}\right) _{\left\vert W_{\xi _{\ast }}\right. }=0$
but $\widetilde{G}\left\vert _{W_{\eta _{\ast }}}\right. \neq 0$, either $%
\mathcal{X}$ is a domain in an algebraic curve of $\mathbb{CP}_{2}$ and has
boundary $\gamma $ in the sense of currents, either $\overline{\mathcal{X}}$
itself is an algebraic curve of $\mathbb{CP}_{2}$ where $\gamma $ is a slit.

\textbf{b2) }If in addition (\ref{F/ sans Green1}) is satisfied by $\Theta
_{\ell }$ for all $c\in H_{1}\left( \mathcal{X}_{\func{reg}}\right) $, then $%
u_{\ell }$ has a (unique) weakly harmonic extension $\widetilde{u_{\ell }}$
and $\Theta _{\ell }=\partial \widetilde{u_{\ell }}$. If $\Theta _{\ell }$
also satisfy $\int_{\mathcal{X}_{\func{reg}}}i\,\Theta _{\ell }\wedge 
\overline{\Theta _{\ell }}<+\infty $, then $\widetilde{u_{\ell }}$ is
harmonic.
\end{theorem}

\noindent \textbf{Remark. }It is possible that $\mathcal{X}$ has zero
boundary in the sense of currents. This occurs only in the exceptional case
where $\overline{\mathcal{X}}$ is a compact complex curve of $\mathbb{CP}%
_{2} $ and (so is algebraic) where $\gamma $ is a slit. In the other cases, $%
\mathcal{X}$ has boundary $\pm \gamma $ in the sense of currents and a
result of Chirka~\cite{ChE1982} gives that outside a zero one Hausdorff
dimensional subset, $\left( \,\overline{\mathcal{X}},\pm \gamma \,\right) $
is locally a manifold with boundary.\medskip

The proofs of the preceding theorems are given in sections three to five.
They use the results on the complex Plateau problem started in \cite%
{WeJ1958c}\cite{BiE1963}, developed in \cite{HaR-LaB1975}\cite{HaR1977}\cite%
{DiT1998b} for $\mathbb{C}^{n}$ and in \cite{HeG1995}\cite{DoP-HeG1997}\cite%
{HaR-LaB2004} for $\mathbb{CP}_{n}$.\medskip

The non constructive existence criterions of theorems~\ref{T/ caract - a}, %
\ref{T/ caract - b} and \ref{T/ caract - c} may incite to seek a less
general but more effective characterization. It has been already mentioned
after theorem~\ref{T/ caract - a} that in the special case $p=0$, the
condition $\left( \partial ^{2}G/\partial \xi _{0}^{2}\right) _{\left\vert
W_{\xi _{\ast }}\right. }=0$ together with the Alexander-Wermer moment
criterion gives an effective tool but only when special DN-data are at hand.

For $p>0$, the main result of \cite{DoP-HeG1997} is that conditions of type (%
\ref{F/ caract G}) and (\ref{F/ SW}) characterize the fact that a given
closed smooth and orientable real chain $\gamma $ of $\mathbb{CP}_{2}$ is,
with adequate orientation, the boundary of some holomorphic chain of $%
\mathbb{CP}_{2}\backslash \gamma $. These conditions have been qualified as
mysterious in \cite{HaR-LaB2004} because the functions satisfying these
relations are produced "deus ex machina". The following criterion, which
completes for a closed connected curve $\gamma $ the one of \cite%
{DoP-HeG1997}, is obtained in \cite{HaR-LaB2004}~: Suppose that the second
coordinate $f_{2}$ of $\mathbb{C}^{2}$ does not vanish on $\gamma $, then
there exists in $\mathbb{CP}_{2}\backslash \gamma $ a connected complex
curve $\mathcal{X}$ with boundary $\pm \gamma $ in the sense of currents if
and only if there exist $p\in \mathbb{N}$ and $A_{d}$ in the space $\mathcal{%
O}\left( d\right) $ of holomorphic homogeneous polynomials of order $d$, $%
1\leqslant d\leqslant p$, such that for $\xi _{0}$ in some neighborhood of $%
0 $, $C_{m}\left( \xi _{0}\right) =\frac{1}{2\pi i}\int_{\gamma }\frac{%
f_{1}^{m}}{f_{2}+\xi _{0}}df_{2}$ satisfies 
\begin{eqnarray*}
C_{k} &=&Q_{k,p}\left( C_{1},...,C_{p}\right) ~~mod~~\mathcal{O}\left(
k\right) ~,~~k>p \\
C_{d}\left( \xi _{0}\right) &=&\dsum\limits_{k>d}\frac{\left( -\xi
_{0}\right) ^{k}}{2\pi i}\int_{\gamma }\frac{f_{1}^{d}}{f_{2}^{k+1}}%
df_{2}+A_{d}\left( \xi _{0}\right) ,\ 1\leqslant d\leqslant p
\end{eqnarray*}%
where $Q_{k,p}$ are universal homogeneous polynomials.

In section~\ref{S/ EffCha}, theorem~\ref{T/ caract - a} is develop for $p>0$
into theorem~\ref{T/ caract - a - Eff} below which gives a more effective
criterion for the Plateau problem in $\mathbb{CP}_{2}$ and also for the
IDN-problem. This new criterion follows from considerations on sums of shock
wave functions modulo affine functions in $\xi _{0}$. Even if decompositions
in sum of shock wave functions are studied for $\xi _{0}$-affine functions,
theorem~\ref{T/ caract - a - Eff} does not consider the case where $G$ is of
that type since it corresponds to a plain case of (\ref{F/ caract G}).

Notations for theorem~\ref{T/ caract - a - Eff}$~$: if $H$ and $u$ are
holomorphic functions on a simply connected domain $D$, we set $\mathcal{D}%
_{H}=\frac{\partial }{\partial \xi _{1}}-\frac{\partial H}{\partial \xi _{0}}
$ and denote by $\mathcal{L}_{H}u$ the unique function $v\in \mathcal{O}%
\left( D\right) $ such that $\partial v/\partial \xi _{0}=\mathcal{D}_{H}u$
and $v\left( 0,.\right) =0$~; $\pi _{1}$ is the projection $\left( \xi
_{0},\xi _{1}\right) \mapsto \xi _{1}$.

%TCIMACRO{\TeXButton{Reprise}{\renewcommand{\theX}{\arabic{X}}}}%
%BeginExpansion
\renewcommand{\theX}{\arabic{X}}%
%EndExpansion

\begin{theorem}
\label{T/ caract - a - Eff}Let $f$\ be defined by (\ref{F/ f}) and consider
the function $G$ defined by (\ref{F/ fonction G}). We assume that $\gamma $
is connected and that $f_{2}$ does not vanish on $\gamma $ so that $G$,
which is assumed to be not affine in $\xi _{0}$, is defined in a simply
connected neighborhood $D$ of $0$ in $\mathbb{C}^{2}.$

\textbf{A/ }If $\left( \overline{\mathcal{X}},\gamma \right) $ is a Riemann
surface with almost smooth boundary where $f$ extends meromorphically, then
the following assertions hold for $G$

\begin{enumerate}
\item There is $p\in \mathbb{N}^{\ast }$ and holomorphic functions $%
a,b,\lambda _{1},...,\lambda _{p-1}$ on $\Delta =\pi _{1}\left( D\right) $
such that the integro-differential equation%
\begin{equation*}
-\mathcal{D}_{G+L}\mathcal{L}_{G+L}^{p-1}\left( G+L\right)
+\dsum\limits_{1\leqslant j\leqslant p-1}\mathcal{D}_{G+L}\mathcal{L}%
_{G+L}^{p-1-j}\,\widetilde{\lambda }_{j}=0,
\end{equation*}%
is valid with $L=\xi _{0}\otimes a+1\otimes b$ and $\widetilde{\lambda }%
_{j}=1\otimes \lambda _{j}$, $1\leqslant j\leqslant p$.

\item For $s_{k}=-\mathcal{L}_{G+L}^{k-1}G+\mathcal{L}_{G+L}^{k-2}\widetilde{%
\lambda }_{1}+\cdots +\mathcal{L}_{G+L}^{0}\widetilde{\lambda }_{k-1}$, $%
1\leqslant k\leqslant p$, the discriminant of $T_{\xi _{0},\xi _{1}}=X^{p}+%
\underset{1\leqslant k\leqslant p}{\Sigma }s_{k}\left( \xi _{0},\xi
_{1}\right) X^{p-k}$ does not vanish identically in $D$.

\item $G=-s_{1}-L$

\item There is $q\in \mathbb{N}$, $\alpha ,\beta \in \mathbb{C}_{q}\left[
\xi _{1}\right] $ such that $\alpha \left( 0\right) =0$, $\deg \beta <q$ and 
$a=\frac{\alpha ^{\,\prime }}{1-\alpha }$, $b=\frac{\beta }{1-\alpha }$
\end{enumerate}

Moreover, if $p$ is the least integer such that (1-2-3) assertions holds, $%
\left( \gamma ,f\right) $ uniquely determines $\left( a,b,\lambda
_{1},...,\lambda _{p-1}\right) $.\smallskip

\textbf{B/ }Assume (1-2-3) holds for some $p\in \mathbb{N}^{\ast }$. Then,
there exists an open Riemann surface $\mathcal{X}$ such that $\overline{%
\mathcal{X}}=\mathcal{X}\cup \gamma $ is a manifold with almost smooth
boundary where $f$ extends meromorphically. Moreover, (4) holds.
\end{theorem}

\noindent \textbf{Remarks. 1. }Non unicity of $\left( a,b,\lambda
_{1},...,\lambda _{p-1}\right) $ solving (1-2-3) means that $\mathcal{X}$
exists but $p$ is not minimal.

\textbf{2. }It is possible that regardless its orientation, $\gamma $ is the
almost smooth boundary of an open Riemann surface $\mathcal{X}$ where $f$
extends meromorphically. It is the case when $\gamma $ cuts a compact
Riemann surface $\mathcal{Z}$ into two smooth domains and $f$ is the
restriction to $\gamma $ of an analytic map from $\mathcal{Z}$ to $\mathbb{CP%
}_{2}$.

\section{Intrinsic EIT on Riemann surfaces\label{S/ Intrinsic}}

The Inverse Dirichlet-Neuman problem, which goes back to Calderon~\cite%
{CaA1980} and which is called now Electrical-Impedance-Tomography problem,
can be sketch like this~: suppose that a bounded domain $\mathcal{X}$ in $%
\mathbb{R}^{2}$ or $\mathbb{R}^{3}$ is an ohmic conductor which means that
the density of current $j$ it may have is proportional (in isotropic cases)
to the electrical field $e=\nabla U$ where $U$ is an electrical potential.
The scalar function $\sigma $ such that $j=\sigma e$ is then called the
conductivity of $\mathcal{X}$ ; $\rho =1/\sigma $ is the resistivity. When
there is no time dependence and no source or sink of current, the equation $%
\limfunc{div}j=0$ holds and Calderon's problem is then to recover $\sigma $
on the whole of $\mathcal{X}$ from the operator $C^{\infty }\left( \gamma
\right) \ni u\mapsto \left( \sigma \nabla \widetilde{u}\right) _{\gamma }$, $%
\widetilde{u}$ being the unique solution of $\limfunc{div}\left( \sigma
\nabla \widetilde{u}\right) =0$ with boundary value $u$.

In what follows, linking the Calderon problem to the Belishev problem
mentioned in the introduction, we formulate the EIT-problem for a more
general setting than the case of domains in $\mathbb{R}^{n}$. The second
part of this section, despite the fact it is also quite elementary, seems to
be new and underline how complex structure is involved in the dimension two
case.\medskip

\noindent \textbf{General dimension}. Assume that $\mathcal{X}$, an open
oriented bordered manifold of dimension $n$ with boundary $\gamma $, is
given with a volume form $\mu $ and a conductivity $\sigma $ modelled as a
tensor from $T^{\ast }\mathcal{X}$ to $\Lambda ^{n-1}T^{\ast }\mathcal{X}$
(see \cite{SyJ1990}). The gradient associated to $\sigma $ relatively to $%
\mu $ is the differential operator which to any $f\in C^{1}\left( \mathcal{X}%
\right) $ associates the tangent vector field $\nabla _{\mu ,\sigma }f$
characterized by%
\begin{equation*}
\left( \nabla _{\mu ,\sigma }f\right) \,\lrcorner \,\mu =\sigma \left(
df\right)
\end{equation*}%
where $\lrcorner $ is the interior product. When $U\in C^{1}\left( \mathcal{X%
}\right) $ is some given potential, the density of physical current $J$ is
by definition$^{\text{(}}$\footnote{%
That $\nabla _{\mu ,\sigma }U$ truly models the density of current is the
assumption of Ohm's law.}$^{\text{)}}$%
\begin{equation*}
J=\nabla _{\mu ,\sigma }U.
\end{equation*}%
If $\mathcal{X}$ has no source or sink of currents and if $U$ has no time
dependence, the flux of current through the boundary of any domain is zero.
Using Stokes formula, this can be modeled by the simplified Maxwell equation%
\begin{equation}
0=\func{div}_{\mu }J=\func{div}_{\mu }\left( \nabla _{\mu ,\sigma }U\right)
\label{F/ E}
\end{equation}%
where $\func{div}_{\mu }$ is$^{\text{(}}$\footnote{%
If $t$ is a differentiable vector field, $\func{div}_{\mu }t$ is defined by $%
d\left( t\,\lrcorner \,\mu \right) =\left( \func{div}_{\mu }t\right) \mu $.}$%
^{\text{)}}$ the divergence with respect to the volume form $\mu $. Going
back to the definition of gradient and divergence, we see that (\ref{F/ E})
is equivalent to the intrinsic equation formulated in \cite{SyJ1990} for
domains in $\mathbb{R}^{n}$%
\begin{equation}
d~\sigma \left( dU\right) =0  \label{F/ Ei}
\end{equation}

Since $\mu $ is no longer involved, the usual DN-operator has to be replaced
by the operator $\Theta $ which to $u\in C^{1}\left( \gamma \right) $
associates $\sigma \left( du\right) _{\gamma }$ which is a section of $%
\Lambda ^{n-1}T_{\gamma }^{\ast }\mathcal{X}$. The Electrical Impedance
Tomography problem, is then to reconstruct $\left( \mathcal{X},\sigma
\right) $ from its DN-map $\Theta $. Of course, the two other aspects of
this problem, unicity and characterization, also has to be studied.

The problem in such a generality is still widely open~; almost all
publications are about domains in $\mathbb{R}^{3}$. In such a case, (\ref{F/
E}) is generally written in euclidean global coordinates. However, when $%
\mathcal{X}$ is a manifold, (\ref{F/ Ei}) yields the same equation in any
chart $\left( W,x\right) $~; setting $\sigma dx_{j}=\underset{1\leqslant
k\leqslant n}{\Sigma }\sigma _{kj}{}\left( -1\right) ^{k}dx_{\widehat{k}}$
with $dx_{\widehat{k}}=\underset{j\neq k}{\wedge }dx_{j}$, (\ref{F/ Ei})
becomes\vspace{-3pt} 
\begin{equation}
\dsum\limits_{1\leqslant k\leqslant n}\dsum\limits_{1\leqslant j\leqslant n}%
\frac{\partial }{\partial x_{k}}\left( \sigma _{k,j}\frac{\partial U}{%
\partial x_{j}}\right) =0  \label{F/ Ec}
\end{equation}

When the conductivity $\sigma $ is symmetric ($\sigma \left( a\right) \wedge
b\equiv \sigma \left( b\right) \wedge a$) and invertible tensor, it is
possible to design a natural metric $g_{\mu ,\rho }$ associated to the
resistivity map $\rho =\sigma ^{-1}$ by the well defined quotient of $n$%
-forms~:%
\begin{equation}
g_{\mu ,\sigma ^{-1}}\left( t\right) =\frac{\sigma ^{-1}\left( t\,\lrcorner
\,\mu \right) \wedge \left( t\,\lrcorner \,\mu \right) }{\mu },~t\in T%
\mathcal{X}\text{.}  \label{F/ Metrique}
\end{equation}%
If $\left( W,x\right) $ is any coordinates chart for $\mathcal{X}$, a direct
calculus in $x$-coordinates shows that for $t=\Sigma t_{k}\partial /\partial
x_{k}$ (\ref{F/ Metrique}) becomes 
\begin{equation}
g_{\mu ,\rho }\left( t\right) =\dsum\limits_{k,\ell }t_{\ell
}\,t_{k}\,\lambda \rho _{k,\ell }  \label{F/ MetParRho}
\end{equation}%
where $\left( \rho _{k,\ell }\right) $ is the matrix of the resistivity $%
\rho =\sigma ^{-1}$ when at any given point $z$ the chosen basis for $%
\Lambda ^{n-1}T_{z}^{\ast }\mathcal{X}$ and $T_{z}^{\ast }\mathcal{X}$ are $%
(\,\left( -1\right) ^{k}dx_{\widehat{k}}\,)$ and $\left( dx_{k}\right) $
respectively. When $\left( \sigma _{j,k}\right) $ is positive definite, $%
g_{\mu ,\rho }$ is a metric on $\mathcal{X}$.

When $n\geqslant 3$, there is a specially adequate choice of metric and
volume.

\begin{proposition}
\label{P/ Hodge}Assume $n\geqslant 3$. Then one can correctly design a
global volume form $\mu $ by letting it be defined by $\mu =\left[ \det
\left( \rho _{k,\ell }\right) \right] ^{\frac{-1}{n-2}}dx_{1}\wedge \cdots
\wedge dx_{n}$ in any coordinates chart $\left( W,x\right) $ for $\mathcal{X}
$. For this specific volume form, $\sigma $ is the Hodge star operator of $%
g_{\mu ,\rho }$ and $\mu $ is the riemannian volume form of $g_{\mu ,\rho }$.
\end{proposition}

This statement, already pointed out by Bossavit and Lee-Uhlmann (see~\cite%
{BoA2001} and~\cite{LaM-UhG2001}) for domains in affine spaces, follows from
calculus in coordinates.

The interest of proposition~\ref{P/ Hodge} is to state the strict
equivalence between the IDN-problem for riemannian manifolds and the
EIT-problem when $n\geqslant 3$. When $\dim \mathcal{X}\geqslant 3$ and $%
\overline{\mathcal{X}}$ is a riemannian real analytic manifold with
boundary, Lassas and Uhlmann have proved in \cite{LaM-UhG2001} that the
DN-operator uniquely determine $\mathcal{X}$ and its metric.\medskip

\noindent \textbf{The two dimensional case}. We now assume $n=2$ and $\sigma
=\rho ^{-1}$ is symmetric and positive so that $\left( \mathcal{X},g_{\mu
,\rho }\right) $ becomes a riemannian manifold whose volume form is
thereafter denoted by $V_{\mu ,\rho }$. Let us emphasize the complex
structure associated to the conformal class of $\left( \mathcal{X},g_{\mu
,\rho }\right) $ by choosing isothermal coordinates charts, that is
holomorphic charts (see e.g. \cite{VeI1962Livre}). In such a chart $\left(
W,z\right) $,%
\begin{equation*}
g_{\mu ,\rho }=\kappa _{\mu ,\rho }\left( dx\otimes dx+dy\otimes dy\right) =%
\func{Re}\left( \kappa _{\mu ,\rho }dz\otimes d\overline{z}\right)
\end{equation*}%
where $x=\func{Re}z$, $y=\limfunc{Im}z$ and $\kappa _{\mu ,\rho }\in
C^{1}\left( W,\mathbb{R}_{+}^{\ast }\right) $. Hence, in these coordinates, $%
\left( \sigma _{k,\ell }\right) =s\func{diag}\left( 1,1\right) $ with $%
\kappa _{\mu ,\rho }=\lambda /s$ and $\lambda \in C^{1}\left( W,\mathbb{R}%
_{+}^{\ast }\right) $ is defined by $\mu =\lambda \,dx\wedge dy=\lambda \,%
\frac{i}{2}dz\wedge d\overline{z}$. Note that $s$ is a global positive
function on $\mathcal{X}$ since it is the well defined quotient of volume
forms~:%
\begin{equation*}
s=\mu /V_{\mu ,\rho }
\end{equation*}%
Note also that $s$ does not depend on $\mu $ and that (\ref{F/ Ei}) finally
evolutes into%
\begin{equation}
d\left( sd^{c}U\right) =0  \label{F/ Elec}
\end{equation}%
where $d^{c}=i\left( \,\overline{\partial }-\partial \,\right) $, $\overline{%
\partial }$ and $\partial $ being the usual global differential operators
associated to the complex structure of the conformal class of $\left( 
\mathcal{X},g_{\mu ,\rho }\right) $. Hence, we have proved the following
which generalizes a result written by Sylvester~\cite{SyJ1990} for domains
in $\mathbb{R}^{2}$.

\begin{proposition}
Let $\mathcal{X}$ be a real two dimensional manifold equipped with a
symmetric and positive tensor $\sigma :T^{\ast }\mathcal{X}\rightarrow
T^{\ast }\mathcal{X}$. Then, there is a complex structure on $\mathcal{X}$
and $s\in C^{1}\left( \mathcal{X},\mathbb{R}_{+}^{\ast }\right) $, called
scalar conductivity, such that (\ref{F/ Ei}) is equivalent to (\ref{F/ Elec}%
).
\end{proposition}

The beginning of this paper has shown that the data $\partial U/\partial \nu 
$ is equivalent to the data $\left( \partial U\right) _{\gamma }$ which
don't involve any metric. Since the knowledge of $\left( \partial U\right)
_{\gamma }$ is equivalent to the knowledge of $\left( sd^{c}U\right)
_{\gamma }$, we consider $\left( sd^{c}U\right) _{\gamma }$ as the DN-datum.
We can now state an intrinsic IDN-problem for two dimensional ohmic
conductors~; for the sake of simplicity, we limit ourselves to manifolds
with boundary and smooth datas.\smallskip

A\textit{\ two dimensional ohmic conductor} is a couple $\left( \mathcal{X}%
,\rho \right) $ where $\mathcal{X}$ is an open oriented bordered two
dimensional real surface (with boundary $\gamma $),the conductivity $\sigma
=\rho ^{-1}$ is a positive definite tensor from $T^{\ast }\mathcal{X}$ to $%
T^{\ast }\mathcal{X}$ and $\mathcal{X}$ is equipped with the complex
structure associated to the riemannian metric $g_{\mu ,\rho }$ defined by (%
\ref{F/ Metrique}) where $\mu $ is any volume form of $\mathcal{X}$. In this
setting, the scalar conductivity is the function $s=\mu /V_{\mu ,\rho }$
where $V_{\mu ,\rho }$ is the volume associated to $g_{\mu ,\rho }$. The
DN-operator is the operator $\theta _{\mathcal{X},\rho }$ defined by%
\begin{equation*}
\theta _{\mathcal{X},\rho }:C^{1}\left( \gamma \right) \ni u\mapsto \left(
sd^{c}\widetilde{u}\right) _{\gamma }\in T_{\gamma }^{\ast }\mathcal{X}
\end{equation*}%
where $\widetilde{u}$ is the unique solution of the following Dirichlet
problem%
\begin{equation}
U\left\vert _{\gamma }\right. =u\hspace{0.5cm}\&\hspace{0.5cm}d\left(
s\,d^{c}U\right) =0.  \label{F/ DiPbDim2}
\end{equation}

The IDN-problem associated to this setting is threefold~:\smallskip

\noindent \textit{Unicity}. Assume that two dimensional ohmic conductors $%
\left( \mathcal{X},\rho \right) $ and $\left( \mathcal{X}^{\prime },\rho
^{\prime }\right) $ share the same boundary $\gamma $ and the same
DN-operator $\theta $. Is it true that there is a diffeomorphism $\varphi :%
\overline{\mathcal{X}}\rightarrow \overline{\mathcal{X}^{\prime }}$ between
manifolds with boundaries such that $\varphi :\mathcal{X}\rightarrow 
\mathcal{X}^{\prime }$ is analytic and $s=s^{\prime }\circ \varphi $ where $%
s $ and $s^{\prime }$ are scalar conductivities of $\mathcal{X}$ and $%
\mathcal{X}^{\prime }$~?\smallskip

\noindent \textit{Reconstruction}. Assume that $\left( \mathcal{X},\rho
\right) $ is a two dimensional ohmic conductor. How from its DN-operator one
can reconstruct a two dimensional ohmic conductor $\left( \mathcal{X}%
^{\prime },\rho ^{\prime }\right) $ which is isomorphic (in the above sense)
to $\left( \mathcal{X},\rho \right) $.\smallskip

\noindent \textit{Characterization}. Let $\gamma $ be a smooth abstract real
curve, $L$ a complex line bundle along $\gamma $ and $\theta $ an operator
from $C^{1}\left( \gamma \right) $ to the space of smooth sections of $L$.
Find a non trivial necessary and sufficient condition on $\left( \gamma
,L,\theta \right) $ which ensures that there exists a two dimensional ohmic
conductor $\left( \mathcal{X},\rho \right) $ such that $L=\Lambda
^{1,0}T_{\gamma }^{\ast }\mathcal{X}$ and $\theta =\theta _{\mathcal{X},\rho
}$.\smallskip

All these problems are open. In the particular case of constant scalar
conductivity $\sigma $, the Dirichlet problem~(\ref{F/ DiPbDim2}) becomes 
\begin{equation*}
U\left\vert _{\gamma }\right. =u~~\&~~dd^{c}U=0
\end{equation*}%
where $dd^{c}=i\partial \overline{\partial }$ is the usual Laplacian. Hence,
our article gives with theorems~\ref{T/ unicite} to~\ref{T/ caract - c} a
rather complete answer to the EIT-problem with constant scalar conductivity.

Concerning the main results given in the literature about unicity,
reconstruction and stability for the important case where $\mathcal{X}$ is a
domain in $\mathbb{R}^{2}$ but the scalar conductivity is not constant, see 
\cite{BrR-UhG1997}\cite{MaN2001} and references therein. Note that the exact
method of reconstruction for this case goes back to \cite{NoR1988}.

\section{Unicity under existence assumption\label{S/ UsousE}}

The notations and hypothesis are taken from theorem and section~\ref{T/
unicite}~; we equip $\mathcal{X}$ with a hermitian metric$^{\text{(}}$%
\footnote{%
Harmonicity does not depend of the chosen hermitian metric.}$^{\text{)}}$ $g$%
. Hence, there is a compact subset $\sigma $ of $\gamma $ such that $%
h^{1}\left( \sigma \right) =0$ and $\left( \overline{\mathcal{X}},\gamma
\right) $ is a manifold with boundary near each point of $\gamma \backslash
\sigma $.

When $u\in C^{\infty }\left( \gamma \right) $, prop.~\ref{L/ Di=Ri} implies
that $u$ has a continuous harmonic extension $\widetilde{u}$ with finite
Dirichlet integral on $\mathcal{X}$. An elementary calculus gives then that
for a fixed continuous defining function $\rho $ of $\gamma $, smooth on $%
\overline{\mathcal{X}}\backslash \sigma $, the operator $L$ defined by (\ref%
{F/ N vers CR et theta}) determines for all $u\in C^{\infty }\left( \gamma
\right) $ the trace on $\gamma $ of the holomorphic $\left( 1,0\right) $%
-form $\partial \widetilde{u}$~: $\partial \widetilde{u}=\left( Lu\right)
\left\vert \partial \rho \right\vert _{g}^{-1}\partial \rho $ on $\gamma
\backslash \sigma $. With (\ref{F/ f}), this implies that $f$ is the
restriction to $\gamma $ of a function $F=\left( F_{1},F_{2}\right) $
meromorphic on $\mathcal{X}$, smooth on $\gamma \backslash \sigma $. Since $%
\left( \gamma ,u,\theta u\right) $ is assumed to be a restricted DN-datum
for $\mathcal{X}$, $F$ is continuous in a neighborhood of $\gamma $ in $%
\overline{\mathcal{X}}$.

The proof of theorem~\ref{T/ unicite} relies on the following lemmas which
enable to see $\overline{\mathcal{X}}$ as a normalization of $\overline{%
\mathcal{Y}}$.

\begin{lemma}
\label{L/ normalisation1}Set $\delta =f\left( \gamma \right) $. Then $%
\mathcal{Y}=F\left( \mathcal{X}\right) \backslash \delta $ is a complex
curve of $\mathbb{C}\mathbb{P}_{2}\backslash \delta $ without compact
component, which has finite mass and satisfies $d\left[ \mathcal{Y}\right] =%
\left[ \delta \right] $. Moreover, each regular point of $\overline{\mathcal{%
X}}$ has in $\overline{\mathcal{X}}$ a neighborhood $V$ such that $%
F:V\rightarrow F\left( V\right) $ is diffeomorphism between manifolds with
smooth boundary.
\end{lemma}

\begin{proof}
Since $F$ is continuous in a neighborhood of $\gamma $ in $\overline{%
\mathcal{X}}$, $\mathcal{Y}$ is a closed set of $\mathbb{CP}_{2}\backslash
\delta $. As $\mathcal{Y}$ is also locally the image of a Riemann surface by
an analytic map, \QTR{EuScript}{$\mathcal{Y}$} is a complex curve of $%
\mathbb{CP}_{2}\backslash \delta $. Since $F_{\ast }\left[ \QTR{EuScript}{%
\mathcal{X}}\right] $ is a locally flat current, the Federer support theorem
(see \cite[p. 316]{HaR1977}\cite[4.1.15 \& 4.1.20]{FeH1969Li}) produces a
locally integrable function $\lambda $ on \QTR{EuScript}{$\mathcal{Y}$} such
that $F_{\ast }\left[ \QTR{EuScript}{\mathcal{X}}\right] =\lambda \left[ 
\QTR{EuScript}{\mathcal{Y}}\right] $ on the regular part \QTR{EuScript}{$%
\mathcal{Y}_{\limfunc{reg}}$ of $\mathcal{Y}$}~; since $d^{2}=0$, $\lambda $
is locally constant. Since $f$ embeds $\gamma $ into $\mathbb{C}^{2}$, each
point $x$ in $\gamma $ which is a regular boundary point of $\overline{%
\mathcal{X}}$ has a neighborhood $V$ of $x$ in $\overline{\mathcal{X}}$ such
that $F:V_{x}\rightarrow F\left( V_{x}\right) $ is a diffeomorphism between
classical manifolds with boundary. Let $F_{W}$ be the restriction of $F$ to
the Riemann surfaces $W=\cup V_{x}$ and $W^{\prime }=F\left( W\right) $. The
degree of $F_{W}$ is at most $1$ otherwise almost all points of $W^{\prime }$
would have at least two different preimages which would imply that $df$ is
zero at almost all points of $\gamma $. So this degree is $1$ and $F_{\ast }%
\left[ \QTR{EuScript}{\mathcal{X}}\right] =\left[ \mathcal{Y}\right] $ on $%
W^{\prime }$. Hence, $\lambda =1$ on each connected component of $\mathcal{Y}
$ and $d\left[ \QTR{EuScript}{\mathcal{Y}}\right] =\left[ \delta \right] $.

If $\mathcal{Y}$ contains a compact complex curve $\mathcal{Z}$, $%
F^{-1}\left( \mathcal{Z}\right) $ is a complex curve in $\mathcal{X}$
without boundary and so is empty. The fact that $\mathcal{Y}$ has a finite
mass follows from a theorem of Wirtinger (see \cite[Lemma 1.5 p. 315]%
{HaR1977}).
\end{proof}

As $\delta $ is smooth, the conclusion of lemma~\ref{L/ normalisation1}
implies, thanks to \cite{HaR-LaB1975}, that $\delta $ contains a compact set 
$\tau $ such that $h^{1}\left( \tau \right) =0$ and $\left( \,\overline{%
\mathcal{Y}},\delta \,\right) $ is manifold with boundary near points of $%
\delta \backslash \tau $. The lemma below described how $\mathcal{Y}$ is
near a point $y$ of $\tau $.

\begin{lemma}
\label{L/ reg bord}Assume $\mathcal{Y}$ is a complex curve of $\mathbb{CP}%
_{2}\backslash \delta $ with finite mass satisfying $d\left[ \mathcal{Y}%
\right] =\pm \left[ \delta \right] $. Let $y$ be a point of $\sigma $ and $U$
a domain containing $y$. Then, among the components of $\mathcal{Y}\cap U$, $%
\mathcal{C}_{y,1}^{U}$, $...$, $\mathcal{C}_{y,m_{U}}^{U}$, one, says $%
\mathcal{C}_{y,1}^{U}$, satisfies $d\left[ \mathcal{C}_{y,1}^{U}\right] =\pm %
\left[ \delta \right] \left\vert _{U_{y}}\right. $ whereas for $j\geqslant 2$%
, $\overline{\mathcal{C}_{y,j}^{U}}\cap U$ is a complex curve of $U$.
\end{lemma}

\begin{proof}
\cite[th. 4.7]{HaR-LaB1975} implies that for each $j$ there is $n_{j}\in 
\mathbb{Z}$ such that $d\left[ \mathcal{\mathcal{C}}_{y,j}^{U}\right] =n_{j}%
\left[ \delta \right] $ on $U$. As $d\left[ \mathcal{Y}\right] =\pm \left[
\delta \right] $, $\Sigma d\left[ \mathcal{C}_{y,j}^{U}\right] =\pm 1$ and
at least one $\mathcal{C}_{y,j}^{U}$, says $\mathcal{C}_{y,1}^{U}$, is such
that $n_{j}\neq 0$. Because $h^{1}\left( \sigma \right) =0$, $\delta \cap U$
contains a point $q$ not in $\sigma $. Then, if $V$ is a sufficiently small
ball centered at $q$, $\mathcal{Y}\cap V$ is submanifold of $V$ with
boundary $\delta \cap V$ and $\mathcal{Y}\cap V$ has only one connected
component which can be nothing else than $\mathcal{C}_{y,1}^{U}\cap V$.
Hence $n_{1}=\pm 1$. Since two different bordered Riemann surfaces of some
open set of $\mathbb{CP}_{2}$ meet at most in a set of zero one dimensional
Hausdorff measure, this implies that $n_{j}=0$ for $j\neq 1$. Thus, if $%
j\geqslant 2$, $d\left[ \mathcal{C}_{y,j}^{U}\right] =0$ and with \cite[th.
2.1 p. 337]{HaR1977} we conclude that $\overline{\mathcal{C}_{y,j}^{U}}\cap
U $ is a complex curve of $U_{y}$.
\end{proof}

If $y\in \delta $, we denote by $m_{y}$ the limit of $m_{U}$ (see lemma~\ref%
{L/ reg bord}) when the diameter of $U$ goes to $0$, $U$ neighborhood of $y$%
~; if $m_{y}\geqslant 2$, then $y\in \tau $. A point $y$ of $\delta $ is
called a \textit{strong singularity of }$\overline{\mathcal{Y}}$ if $y$ is
not a regular point of $\mathcal{C}_{y,1}^{U}$ and a \textit{weak singularity%
} of $\overline{\mathcal{Y}}$ if $y$ is regular point of $\mathcal{C}%
_{y,1}^{U}$ but $m_{y}\geqslant 2$.

We denote by $\tau _{1}$ (resp. $\tau _{2}$) the sets of points where $%
\overline{\mathcal{Y}}$ has weak (resp. strong) singularity. Then $\tau
=\tau _{1}\cup \tau _{2}$ and $\tau _{1}\cap \tau _{2}=\varnothing $. Note
that both $\tau _{1}$ and $\tau _{2}$ may contains points $y$ where $%
m_{y}\geqslant 2$.

We denote by $\overline{\mathcal{Y}}_{\func{sing}}=\mathcal{Y}_{\func{sing}%
}\cup \tau $ the singular locus of $\overline{\mathcal{Y}}$, that is the set
of points of $\overline{\mathcal{Y}}$ where $\overline{\mathcal{Y}}$ is not
a smooth manifold with boundary and we set $\mathcal{B}=f\left( \sigma
\right) \cup \overline{\mathcal{Y}}_{\func{sing}}$, $\mathcal{A}%
=F^{-1}\left( \mathcal{B}\right) $ and $\mathcal{X}_{\circ }=\mathcal{X}%
\backslash F^{-1}\left( \delta \right) =\mathcal{X}\backslash F^{-1}\left(
\tau _{2}\right) $.

\begin{lemma}
\label{L/ normalisation2}The map $F:\overline{\mathcal{X}}\rightarrow 
\overline{\mathcal{Y}}$ is a normalization in the following sense~: $F:%
\mathcal{X}_{\circ }\rightarrow \mathcal{Y}$ is a (usual) normalization and $%
F:\overline{\mathcal{X}}\backslash \mathcal{A}\rightarrow \overline{\mathcal{%
Y}}\backslash \mathcal{B}$ is a diffeomorphism between manifolds with
boundary.
\end{lemma}

\begin{proof}
Since $\mathcal{X}_{\circ }=\mathcal{X}\backslash F^{-1}\left( \delta
\right) $, the properness of $F\left\vert _{\mathcal{X}_{\circ }}\right. $
and the finiteness of its fibers are elementary. For each connected
component $\mathcal{C}$ of $\mathcal{X}_{\circ }\backslash \mathcal{A}$, the
degree $m_{\mathcal{C}}$ of $F:\mathcal{C}\rightarrow F\left( \mathcal{C}%
\right) $ as a Riemann surfaces morphism is finite and $F_{\ast }\left[ 
\mathcal{C}\right] =\delta _{\mathcal{C}}\left[ F\left( \mathcal{C}\right) %
\right] $. Reasoning as in lemma~\ref{L/ normalisation1}'s proof, we get $m_{%
\mathcal{C}}=1$. As $\mathcal{Y}_{\func{sing}}$ contains all the points of $%
\mathcal{Y}$ which has more than one preimage by $F$, $F:\mathcal{C}%
\rightarrow F\left( \mathcal{C}\right) $ is thus an isomorphism. Let $%
\mathcal{C}^{\prime }$ be another connected component of $\mathcal{X}_{\circ
}\backslash \mathcal{A}$ and assume that $F\left( \mathcal{C}\right) $ and $%
F\left( \mathcal{C}^{\prime }\right) $ meet at $q$. Since $q\notin \mathcal{B%
}$, the germs of $F\left( \mathcal{C}\right) $ and $F\left( \mathcal{C}%
^{\prime }\right) $ at $q$ are equal. This leads to $F\left( \mathcal{C}%
\right) =F\left( \mathcal{C}^{\prime }\right) $ which yields the
contradiction $d\left[ \mathcal{Y}\right] =2\left[ \delta \right] $ near
regular boundary points of $\delta $ in $bF\left( \mathcal{C}\right) $.
Hence, $F:\mathcal{X}_{\circ }\backslash \mathcal{A}\rightarrow \mathcal{Y}%
\backslash \mathcal{B}=\mathcal{Y}_{\func{reg}}$ is an isomorphism of
complex manifolds. As $\mathcal{X}_{\circ }\cap \mathcal{A}=F^{-1}\left( 
\mathcal{Y}_{\func{sing}}\right) $ has empty interior, $F:\mathcal{X}_{\circ
}\rightarrow \mathcal{Y}$ is a usual normalization.

Set $\widetilde{\mathcal{X}}=\overline{\mathcal{X}}\backslash \mathcal{A}$, $%
\widetilde{\mathcal{Y}}=\overline{\mathcal{Y}}\backslash \mathcal{B}$, $%
\widetilde{\tau }=\tau \cup f\left( \sigma \right) $ and $\widetilde{\delta }%
=\delta \backslash \widetilde{\tau }$~; by definition of $\mathcal{B}$, $%
\widetilde{\mathcal{Y}}$ is a manifold with smooth boundary $\widetilde{%
\delta }=\delta \backslash \widetilde{\tau }$ and $\widetilde{\mathcal{X}}$
has smooth boundary $\gamma \backslash \widetilde{\sigma }$ where $%
\widetilde{\sigma }=f^{-1}\left( \tau \right) =\sigma \cup f^{-1}\left( \tau
\right) $. The map $F:\widetilde{\mathcal{X}}\rightarrow \widetilde{\mathcal{%
Y}}$ is onto by construction. It is injective because the maps $F:\mathcal{X}%
_{\circ }\rightarrow \mathcal{Y}$ and $f:\gamma \rightarrow \delta $ are so
and because if $x_{1}\in \mathcal{X}$ and $x_{2}\in \gamma $ have the same
image $y$ by $F$, then $m_{y}\geqslant 2$, $y\in \tau $ and $x_{1},x_{2}\in 
\mathcal{A}$. Since $\widetilde{\mathcal{Y}}\backslash \widetilde{\tau }=%
\mathcal{Y}_{\func{reg}}$ and $F:\mathcal{X}_{\circ }\backslash \mathcal{A}%
\rightarrow \mathcal{Y}_{\func{reg}}$ is a diffeomorphism, the fact that $F:%
\widetilde{\mathcal{X}}\rightarrow \widetilde{\mathcal{Y}}$ is a
diffeomorphism between manifolds with boundary has only to be check locally
near boundary points. If $x\in \gamma \backslash \widetilde{\sigma }$, then $%
y=f\left( x\right) \notin \tau $ and the last conclusion of lemma~\ref{L/
normalisation1} implies that there are open neighborhoods $V$ and $W$ of $x$
and $y$ in $\overline{\mathcal{X}}$ and $\overline{\mathcal{Y}}$ such that $%
F:V\rightarrow W$ is diffeomorphism between manifolds with boundary.
\end{proof}

\subsection{Proof of theorem~\protect\ref{T/ unicite}.}

Let $L^{\prime }$ be the operator defined by (\ref{F/ N vers CR et theta})
when $N$ is changed for $N^{\prime }$, let us denote $F^{\prime }$ the
meromorphic extension of $f$ to \QTR{EuScript}{$\mathcal{X}$}$^{\prime }$
and let $\mathcal{Y}^{\prime }=F^{\prime }\left( \QTR{EuScript}{\mathcal{X}}%
^{\prime }\right) \backslash \delta $ where $\delta =f\left( \gamma \right) $%
. By lemma~\ref{L/ normalisation1}, the sets $\mathcal{Y}^{\prime }$ and $%
\mathcal{Y}$ are two complex curves of $\mathbb{CP}_{2}\backslash \delta $
which has no compact component and both are bordered by $\left[ \delta %
\right] $ in the sense of currents. Hence they are identical by a
consequence of a Harvey-Shiffman theorem (see \cite[prop. 1.4.1]{DoP-HeG1997}%
).

Taking in account lemma~\ref{L/ normalisation2} and the fact that $\mathcal{B%
}\cap \mathcal{Y}=\mathcal{Y}_{\func{sing}}=\mathcal{B}^{\prime }\cap 
\mathcal{Y}$, this implies that $\Phi =F^{-1}\circ F^{\prime }$ is an
analytic isomorphism between \QTR{EuScript}{$\mathcal{X}$}$^{\prime
}\backslash F^{\prime -1}\left( \mathcal{Y}_{\func{sing}}\right) $ and 
\QTR{EuScript}{$\mathcal{X}$}$^{\prime }\backslash F^{-1}\left( \mathcal{Y}_{%
\func{sing}}\right) $. Using the properness of $F:\mathcal{X}_{\circ
}\rightarrow \mathcal{Y}$ and $F^{\prime }:\mathcal{X}_{\circ }^{\prime
}\rightarrow \mathcal{Y}$, we conclude that $\Phi $ extends holomorphically
to \QTR{EuScript}{$\mathcal{X}^{\prime }$}. Likewise, $\Psi =F^{\prime
-1}\circ F$ extends holomorphically to \QTR{EuScript}{$\mathcal{X}$}. As $%
\Phi (\Psi (x^{\prime }))=x^{\prime }$ and $\Psi (\Phi (x))=x$ for almost
all $x^{\prime }\in \mathcal{X}^{\prime }$ and $x\in \mathcal{X}$, the
extension of $\Phi $ is an isomorphism from $\mathcal{X}$ to $\mathcal{X}%
^{\prime }$.

As $F$ and $F^{\prime }$ extend $f$ to $\mathcal{X}$ and $\mathcal{X}%
^{\prime }$, $\Phi $ extend continuously to $\gamma $ by the identity map on 
$\gamma $ . Set $\sigma =\overline{\mathcal{X}}_{\func{sing}}$ and $\sigma
^{\prime }=\overline{\mathcal{X}^{\prime }}_{\func{sing}}$ and let $%
x^{\prime }$ be in $\gamma \backslash \left( \sigma \cup \sigma ^{\prime
}\right) $. Then if $y=f\left( x\right) \notin \tau =\delta \cap \overline{%
\mathcal{Y}}_{\func{sing}}$, $\Phi $ is a diffeomorphism between
neighborhoods of $x$ in $\overline{\mathcal{X}}$ and $\overline{\mathcal{X}%
^{\prime }}$ because $F$ (resp. $F^{\prime }$) is a diffeomorphism between
manifold with boundary from a neighborhood of $x$ in $\overline{\mathcal{X}}$
(resp. in $\overline{\mathcal{X}^{\prime }}$) to a neighborhood of $y$ in $%
\overline{\mathcal{Y}}$. If $y\in \tau $, then the last conclusion of lemma~%
\ref{L/ normalisation1} implies that $y\in \tau _{1}$ so that there is a
open neighborhood $U$ of $y$, a component $\mathcal{C}_{y,1}^{U}$ of $%
\mathcal{Y}\cap U$ and open neighborhoods $V$ and $V^{\prime }$ of $x$ in $%
\overline{\mathcal{X}}$ and $\overline{\mathcal{X}^{\prime }}$ such that
that $F:V\rightarrow \mathcal{C}_{y,1}^{U}$ and $F^{\prime }:V^{\prime
}\rightarrow \mathcal{C}_{y,1}^{U}$ are diffeomorphism between manifolds
with smooth boundary. Hence, $\Phi :V\rightarrow V^{\prime }$ is
diffeomorphism between manifolds with smooth boundary. Finally, $\Phi $
realizes a diffeomorphism between manifolds with smooth boundary from $%
\mathcal{X}^{\prime }\backslash \left( \sigma \cup \sigma ^{\prime }\right) $
to $\mathcal{X}\backslash \left( \sigma \cup \sigma ^{\prime }\right) $ and
the proof is complete.$~~_{\blacksquare }$\medskip

The proof contains the following variation of theorem~\ref{T/ unicite}.

\begin{theorem}
\label{T/ unicite de Y}Assume that $\mathcal{X}$ and $\mathcal{X}^{\prime }$
are open Riemann surfaces with almost smooth boundary $\gamma $ such that
the map $f$ defined by (\ref{F/ f}) is an embedding of $\gamma $ into $%
\mathbb{CP}_{2}$ and has a meromorphic extension $F$ to $\mathcal{X}$ and $%
F^{\prime }$ to $\mathcal{X}^{\prime }$ which are continuous near $\gamma $.
Then $F\left( \mathcal{X}\right) \backslash f\left( \gamma \right)
=F^{\prime }\left( \mathcal{X}^{\prime }\right) /f\left( \gamma \right) 
\overset{def}{=}\mathcal{Y}$ is a complex curve of $\mathbb{C}\mathbb{P}%
_{2}\backslash \delta $ without compact component, which has finite mass and
satisfies $d\left[ \mathcal{Y}\right] =\left[ \delta \right] $. Moreover, $%
\overline{\mathcal{X}}$ and $\overline{\mathcal{X}^{\prime }}$ are
normalizations of $\overline{\mathcal{Y}}$ in the sense of lemma~\ref{L/
normalisation2}.
\end{theorem}

Thus, Riemann surfaces constructed in the converse part of theorems~\ref{T/
caract - a}, \ref{T/ caract - b} and \ref{T/ caract - c} are the only
possible candidates for a solution to the IDN-problem.

\section{Existence and reconstruction, proof of theorem~\protect\ref{T/
reconstr}}

We first prove that the Stokes formula holds in almost smoothly bordered
manifolds; one can see also \cite{FeH1969Li}.

\begin{lemma}
\label{L/ Stokes}Let $\left( \,\overline{\mathcal{X}},\gamma \,\right) $ be
a Riemann surface with almost smooth boundary. Then for any 1-form $\varphi $
which is continuous on $\overline{\mathcal{X}}$ such that $d\varphi $ exists
as an integrable differential on $\mathcal{X}$, we have%
\begin{equation}
\int_{\mathcal{X}}d\varphi =\int_{\gamma }\varphi  \label{F/ Stokes}
\end{equation}
\end{lemma}

\begin{proof}
Set $\sigma =\overline{\mathcal{X}}_{\func{sing}}$. Since $h^{2}\left( \,%
\overline{\mathcal{X}}\,\right) <\infty $ and $h^{1}\left( \sigma \right) =0$%
, there is an increasing sequence $\left( \mathcal{X}_{k}\right) $ of smooth
open sets of $\mathcal{X}$ such that $\left( h^{1}\left( b\mathcal{X}%
_{k}\right) \right) $ and $\left( h^{2}\left( \,\overline{\mathcal{X}}%
\backslash \mathcal{X}_{k}\right) \,\right) $ both have limit zero and $%
\overline{\mathcal{X}}\backslash \overline{\mathcal{X}_{k}}$ is contained in
a $2^{-k}$-neighborhood of $\sigma $. Let $\varphi $ be as above. Since $%
d\varphi $ is integrable and $\lim h^{2}\left( \mathcal{X}_{k}\right) =0$, $%
\left( \int_{\mathcal{X}_{k}}d\varphi \right) $ has limit $\int_{\mathcal{X}%
}d\varphi $. As $b\mathcal{X}_{k}=\left( \gamma \cap \overline{\mathcal{X}%
_{k}}\right) \cup \left[ \left( b\mathcal{X}_{k}\right) \backslash \gamma %
\right] $ and $\lim h^{1}\left( b\mathcal{X}_{k}\right) =0$, $\left(
\int_{\left( b\mathcal{X}_{k}\right) \backslash \gamma }\varphi \right) $
converges to $0$. Hence, $\lim \int_{\gamma \cap \overline{\mathcal{X}_{k}}%
}\varphi =\int_{\gamma }\varphi $ and the classical Stokes formula for $%
\varphi $ and $\mathcal{X}_{k}$ yields (\ref{F/ Stokes}).
\end{proof}

We prove now a variation of the Riemann's existence theorem.

\begin{proposition}
\label{L/ Di=Ri}Let $\left( \,\overline{\mathcal{X}},\gamma \,\right) $ be a
Riemann surface with almost smooth boundary and $u$ a real valued
lipschitzian function on $\gamma $. Then $u$ has a unique continuous
harmonic extension $\widetilde{u}$ of $u$ to $\mathcal{X}$ and $\widetilde{u}
$ has finite Dirichlet integral $\int i\,\partial \widetilde{u}\wedge 
\overline{\partial }\widetilde{u}$. Moreover, $N_{\mathcal{X}}u$ defined as $%
\partial \widetilde{u}/\partial \nu $ on $\gamma \backslash \overline{%
\mathcal{X}}_{\func{sing}}$ admits an extension on $\gamma $ as a current of
order $1$ on $\gamma $.
\end{proposition}

\begin{proof}
Following the lines of Riemann's method for harmonic extension of smooth
functions, we first construct an adequate space $W^{1}\left( \mathcal{X}%
\right) $.

Since $\left( \,\overline{\mathcal{X}},\gamma \,\right) $ is at least a
topological bordered manifold, for every fixed point $x$ in $\gamma $ we can
choose in $\mathcal{X}$ an open set $\Delta _{x}$ whose closure in $%
\overline{\mathcal{X}}$ is a neighborhood of $x$ and which is mapped by a
complex coordinate $\varphi _{x}$ into the closure of the unit disk $\mathbb{%
D}$ of $\mathbb{C}$, $\varphi _{x}$ being a homeomorphism from $\overline{%
\Delta _{x}}$ to $\mathbb{\overline{D}}$. Note that if $x^{\prime }\in
\gamma \cap \overline{\Delta _{x}}$ is a regular point of $\overline{%
\mathcal{X}}$, $\varphi _{x}$ has to be diffeomorphism between manifolds
with boundary from a neighborhood of $x^{\prime }$ to a neighborhood of $%
\varphi _{x}\left( x^{\prime }\right) $ in $\overline{\mathbb{D}}$. If $x\in 
\mathcal{X}$, we choose a conformal open disk $\varphi _{x}:\Delta
_{x}\rightarrow \mathbb{D}$ of $\mathcal{X}$ centered at $x$. With the help
of a continuous partition of unity, we can now construct a continuous
hermitian metric $h$ on $\mathcal{X}$ by gluing together the local metrics $%
\left( \varphi _{x}\right) _{\ast }dz\wedge d\overline{z}$ where $z$ is the
standard coordinate of $\mathbb{C}$. We then denote by $W^{1}\left( \mathcal{%
X}\right) $ the Sobolev space of functions in $L^{2}\left( \mathcal{X}%
,h\right) $ with finite Dirichlet integral.

By construction, any function $A$ in $W^{1}\left( \mathcal{X}\right) $ is
such that for each $x\in \gamma $, $B_{x}=\left( \varphi _{x}\right) _{\ast
}A\left\vert _{\Delta _{x}}\right. $ is square integrable for the standard
metric of $\mathbb{D}$. Since the values of Dirichlet integrals are
conformal invariants, it follows that $B_{x}$ is in the standard Sobolev
space $W^{1}\left( \mathbb{D}\right) $ and hence admits a boundary value $%
b_{x}$ on $\mathbb{T}=b\mathbb{D}$ which is in $W^{1/2}\left( \mathbb{T}%
\right) $. As $b_{x}$ is punctually defined almost everywhere, $%
a_{x}=b_{x}\circ \varphi $ is defined almost everywhere in $\gamma \cap 
\overline{\Delta _{x}}$. The constructions made for each $x\in \gamma $ glue
together to form a function defined almost everywhere in $\gamma $ which we
call the boundary value of $A$.\smallskip

We consider now the subset $F$ of $W^{1}\left( \mathcal{X}\right) $ with
boundary value $u$. It is closed and non empty since by a result of McShane~%
\cite{McE1934}, $u$ admits a Lipschitz extension to $\mathcal{X}$. It
follows now from classical arguments that the Dirichlet integral can be
minimized in $F$ at some function $\widetilde{u}$ which has to be harmonic
in $\mathcal{X}$. It remains only to show that $\widetilde{u}$ is continuous
on $\overline{\mathcal{X}}$. If $x\in \gamma $, what precedes implies that $%
v_{x}=\left( u_{x}\circ \varphi _{x}^{-1}\right) \left\vert _{\mathbb{T}%
}\right. $ is in $W^{1/2}\left( \mathbb{T}\right) $, continuous near $%
\varphi _{x}\left( x\right) $ and is the boundary value of $\widetilde{v}%
_{x}=\widetilde{u}\circ \varphi _{x}^{-1}$. Hence, classical Poisson formula
for the disc implies that near $\varphi _{x}\left( x\right) $ in $\overline{%
\mathbb{D}}$, $\widetilde{v}_{x}$ is continuous up to $\mathbb{T}$ with
restriction $v_{x}$ on $\mathbb{T}$. Since $\varphi _{x}$ is an
homeomorphism, we get that $\widetilde{u}$ is continuous at $x$ with value $%
u\left( x\right) $.

Let $\theta u$ be the form defined by~\ref{F/ N vers CR et theta}. The
Stokes formula (\ref{F/ Stokes}) implies that if $\varphi \in C^{1}\left(
\gamma \right) $ and $\Phi $ is a Lipschitz extension of $\varphi $ on $%
\mathcal{X}$, 
\begin{equation*}
\int_{\gamma }\varphi \theta u=-\int_{\mathcal{X}}\partial \widetilde{u}%
\wedge \overline{\partial }\Phi .
\end{equation*}%
As the last integral is independent of the Lipschitz extension of $\varphi $%
, this means that $\theta u$ and hence $N_{\mathcal{X}}u$, are well defined
currents of order $1$.
\end{proof}

Assume now that the hypothesis of theorem~\ref{T/ reconstr} are true. Lemma~%
\ref{L/ normalisation1} points out that $F$ projects $\gamma $ on a smooth
curve $\delta $ of $\mathbb{C}^{2}$ which bounds in the sense of currents a
complex curve $\mathcal{Y}$ of $\mathbb{CP}_{2}\backslash \delta $ which has
finite mass and no compact component and theorem~\ref{T/ unicite de Y}
implies that for some subset $\mathcal{X}_{\circ }$ of $\mathcal{X}$ with
discrete complement in $\mathcal{X}$, $F:\mathcal{X}_{\circ }\rightarrow 
\mathcal{Y}$ is usual normalization. Hence, if $\mathcal{B}=\mathcal{Y}_{%
\func{sing}}$ and $\mathcal{A}=F^{-1}\left( \mathcal{B}\right) $, $F:%
\mathcal{X}\backslash \mathcal{A}\rightarrow \mathcal{Y}\backslash \mathcal{B%
}$ is one to one. This is part 1 of theorem~\ref{T/ reconstr}.\medskip

Before proving the second claim of theorem~\ref{T/ reconstr}, we recall that 
$\mathbb{CP}_{2}$ is equipped with homogenous coordinates $w$ and $\mathbb{C}%
^{2}$ identified with $\left\{ w_{0}\neq 0\right\} $ has affine coordinates $%
z_{1}=w_{1}/w_{0}$ and $z_{2}=w_{2}/w_{0}$. Set $\Delta _{\infty }=\left\{
w_{0}=0\right\} $, $\mathcal{Y}_{\infty }=\mathcal{Y}\cap \Delta _{\infty }$
and if $\xi \in \mathbb{C}$, we set $\Delta _{\xi }=\{w_{2}=\xi w_{0}\}$ and 
$\mathcal{Y}_{\xi }=\mathcal{Y}\cap \Delta _{\xi }$. Set 
\begin{equation*}
\Omega _{\xi }^{m}=\frac{z_{1}^{m}}{z_{2}-\xi }dz_{2}=\frac{w_{1}^{m}}{%
w_{0}^{m}}\frac{dw_{2}}{w_{2}-\xi w_{0}}-\frac{w_{1}^{m}}{w_{0}^{m+1}}\frac{%
w_{2}dw_{0}}{w_{2}-\xi w_{0}}.
\end{equation*}%
Applying the Stokes formula either for $\mathcal{Y}$ or $\mathcal{X}$, it
turns out that $\frac{1}{2\pi i}\int_{\gamma }\frac{f_{1}^{m}}{f_{2}-\xi }%
df_{2}$ equals $S_{m}\left( \xi \right) +P_{m}\left( \xi \right) $ where 
\begin{equation*}
S_{m}\left( \xi \right) =\dsum\limits_{z\in \mathcal{Y}_{\xi }}\limfunc{Res}%
\left( \eta ^{\ast }\Omega _{\xi }^{m},z\right) ~~;~~~P_{m}\left( \xi
\right) =\dsum\limits_{z\in \mathcal{Y}_{\infty }}\limfunc{Res}\left( \eta
^{\ast }\Omega _{\xi }^{m},z\right)
\end{equation*}%
and $\eta :\mathcal{Y}\rightarrow \mathbb{CP}_{2}$ is the canonical
injection.

For almost all $\xi _{\ast }$ in $\mathbb{C}$, $\mathcal{Y}$ meets
transversely $\Delta _{\xi _{\ast }}$ only in $\mathbb{C}^{2}\cap \mathcal{Y}%
_{\func{reg}}$~; for such a fixed $\xi _{\ast }$, set $p=\func{Card}\mathcal{%
Y}_{\xi _{\ast }}$ and $\mathcal{Y}_{\xi _{\ast }}=\{z_{1\ast },...,z_{p\ast
}\}$. For $\xi $ in a sufficiently small connected neighborhood $W_{\xi
_{\ast }}$ of $\xi _{\ast }$, $\mathcal{Y}_{\xi }$ lies then in $\mathbb{C}%
^{2}\cap \mathcal{Y}_{\func{reg}}$ and can be written $\left\{ z_{1}\left(
\xi \right) ,...,z_{p}\left( \xi \right) \right\} $ with $z_{j}\left( \xi
\right) =\left( h_{j}\left( \xi \right) ,\xi \right) $ where $h_{j}$ is
holomorphic in $W_{\xi _{\ast }}$ and has value $z_{j\ast }$ at $\xi _{\ast
} $, $1\leqslant j\leqslant p$. Direct calculation shows (see \cite%
{DoP-HeG1997}) that the poles of $\Omega _{\xi }^{m}$ in $\mathbb{C}^{2}$
are $z_{1}\left( \xi \right) ,...,z_{p}\left( \xi \right) $ with residue $%
h_{1}\left( \xi \right) ^{m},...,h_{p}\left( \xi \right) ^{m}$. Hence, $%
S_{m}=S_{h,m}$ in $V$.

Reasoning as in lemma~\ref{L/ diffaffine} in next section, we can assume
without loss of generality that $\mathcal{Y}$ meets transversely $\Delta
_{\infty }$ and that $\mathcal{Y}_{\infty }\subset \mathcal{Y}_{\func{reg}}$%
. In this situation, a direct calculus (see \cite{DoP-HeG1997}) gives that
at $y\in \mathcal{Y}_{\infty }$, $\Omega _{\xi }^{m}$ has a pole of order $%
m+1$ with a residue which is a polynomial in $\xi $ of degree at most $m$.
Hence, $P_{m}$ is is a polynomial in $\xi $ of degree at most $m$~; formula (%
\ref{F/ equation pour Y}) is proved.\smallskip

If $A\geqslant B$ and $\xi _{0},...,\xi _{A-1}$ are mutually distinct, the
Vandermonde matrix $\left( \xi _{\nu }^{\mu }\right) _{0\leqslant \nu ,\mu
\leqslant B-1}$ is invertible and hence, the system $(E_{m,\xi _{\nu
}})_{0\leqslant \nu \leqslant B-1}$ enable to write the coefficients of $%
P_{m}$ as a linear combination of the $S_{h,m}\left( \xi _{\nu }\right) $, $%
0\leqslant \nu \leqslant B-1$. Introducing this result in $(E_{\xi
})=(E_{m,\xi _{\nu }})_{\substack{ 0\leqslant m\leqslant B-1  \\ 0\leqslant
\nu \leqslant B-1}}$ we get a linear system which, since $AB-\frac{1}{2}%
B\left( B+1\right) \geqslant pA$ when $B\geqslant 2p+1$, enable to compute
for a generic $\xi $ the unknowns $S_{h,m}\left( \xi _{\nu }\right) $ and,
thanks to the Newton-Girard formulas, the elementary symmetric functions of $%
h_{1}\left( \xi _{\nu }\right) ,...,h_{p}\left( \xi _{\nu }\right) $~;
finally we get the intersection points $\left( h_{j}\left( \xi _{\nu
}\right) ,\xi _{\nu }\right) $ of $\mathcal{Y}$ with $\Delta _{\xi _{\nu }}$%
.\medskip

We prove the third assertion of theorem~\ref{T/ reconstr}. Almost all $\xi
_{\ast }$ in $\mathbb{C}^{2}$ has a connected neighborhood $W_{\xi _{\ast }}$
such that there is a compact of $\mathbb{CP}_{2}\backslash \left( \delta
\cup \mathcal{Y}_{\func{sing}}\right) $ containing all $\mathcal{Y}_{\xi }$
when $\xi \in W_{\xi _{\ast }}$. When $\xi _{\ast }$ is such and $\xi \in
W_{\xi _{\ast }}$ the form $\Phi _{\xi }^{m,\ell }=\frac{F_{1}^{m}}{%
F_{2}-\xi }\partial \widetilde{u_{\ell }}$ may have poles of order at most $%
m $ at infinity i.e. in $\left\{ w_{0}=0\right\} \cap \mathcal{X}$ and
whiles its over poles lies in a compact of $\mathcal{X}$. Since $\widetilde{u%
}_{\ell }$ is the continuous harmonic extension of $u_{\ell }$ on $\mathcal{X%
} $, $\int_{\mathcal{X}}i\,\partial \widetilde{u}_{\ell }\wedge \overline{%
\partial }\widetilde{u}_{\ell }<+\infty $ by proposition~\ref{L/ Di=Ri} and
we can apply the Stokes formula (\ref{F/ Stokes}) to it on $\mathcal{X}$.
This gives (\ref{F/ equation Theta}) after a residue calculus.\medskip

\noindent \textbf{Remark}. The $\frac{1}{2}B\left( B+1\right) $ coefficients
of the polynomials $P_{k}$ come from the residues of the intersection points
of $\mathcal{Y}$ with $\left\{ w_{0}=0\right\} $. In the generic case, $%
\mathcal{Y}$ is given near theses points as the graph of holomorphic
functions $\psi _{1},...,\psi _{q}$ of the variable $w_{0}/w_{2}$, it
appears that the coefficients of $P_{k}$ are ruled by the derivatives of
order at most $k$ at $0$ of the $\psi _{\ell }$. The reconstruction of $%
\mathcal{X}$ is thus possible with a non linear system with only $pA+q\left(
B+1\right) $ unknowns.

\section{Proofs of characterizations theorem~\protect\ref{T/ caract - a}, 
\protect\ref{T/ caract - b} and \protect\ref{T/ caract - c}}

The proofs of theorems~\ref{T/ caract - a}, \ref{T/ caract - b} and~\ref{T/
caract - c} follow a similar schema. The function $f$ defined by (\ref{F/ f}%
) embeds $\gamma $ into a smooth real curve $\delta =f\left( \gamma \right) $
of $\mathbb{C}^{2}$. The necessary conditions for the existence of a
solution to the IDN-problem for $\gamma $ are drawn from the fact this
existence implies that $\delta $ bounds a "concrete" Riemann surface in $%
\mathbb{CP}_{2}$ or $\mathbb{C}^{2}$. The sufficient part of theorem~\ref{T/
caract - c} reconstructs the concrete but singular solution to the
IDN-problem~; a normalization gives then the sufficient part of ~\ref{T/
caract - b}. The proof of theorem$~$\ref{T/ caract - a} follows a similar
scheme.

\subsection{Proof of A.theorem~\protect\ref{T/ caract - a}\label{S/
preuve3aA}}

Assume that $\mathcal{X}$ is an open bordered riemannian surface of finite
volume with restricted DN-datum $\left( \gamma ,u,\theta u\right) $. Then
the functions $F_{j}$ ($j=1,2$) which are the well defined quotients of
forms $\left( \partial \widetilde{u_{j}}\right) /\left( \partial \widetilde{u%
}_{0}\right) $ are meromorphic and letting $F=\left( F_{1},F_{2}\right) $,
lemma~\ref{L/ normalisation1} implies that $\mathcal{Y}=F\left( \mathcal{X}%
\right) \backslash \delta $, $\delta =f\left( \gamma \right) $, is a complex
curve of finite volume, without compact component and bordered by $\left[
\delta \right] $ in the sense of currents. Moreover, the function $G$ has
the expression%
\begin{equation*}
G\left( \xi _{0},\xi _{1}\right) =\frac{1}{2\pi i}\int_{\delta }\Omega _{\xi
}~,~\Omega _{\xi }=\frac{w_{1}}{w_{0}}\frac{d\Lambda _{\xi }\left( w\right) 
}{\Lambda _{\xi }\left( w\right) }-\frac{w_{1}}{w_{0}^{2}}dw_{0}~,
\end{equation*}%
where $\delta =f\left( \gamma \right) $, $\left( w_{0}:w_{1}:w_{2}\right) $
are homogenous coordinates for $\mathbb{CP}_{2}$ and $\Lambda _{\xi }\left(
w\right) =\xi _{0}w_{0}+\xi _{1}w_{1}+w_{2}$.

For almost all $\xi _{\ast }=\left( \xi _{0\ast },\xi _{1\ast }\right) $ and
for all $\xi $ in a sufficiently small connected neighborhood $W_{\xi _{\ast
}}$ of $\xi _{\ast }$, $\mathcal{Y}$ meets $\Delta _{\xi }=\left\{ \Lambda
_{\xi }=0\right\} $ transversely, $\mathcal{Y}_{\xi }=\mathcal{Y}\cap \Delta
_{\xi }\subset \mathbb{C}^{2}\cap \mathcal{Y}_{\func{reg}}$ so that there
exists $p=\func{Card}\mathcal{Y}_{\xi _{\ast }}$ holomorphic functions $%
H_{j}=\left( 1:h_{j}:h_{j,2}\right) :W_{\xi _{\ast }}\rightarrow \mathbb{CP}%
_{2}$ such that $\mathcal{Y}_{\xi }=\left\{ H_{j}\left( \xi \right)
,~1\leqslant j\leqslant p\right\} $ and $\left( h_{j}\right) _{1\leqslant
j\leqslant p}$ are mutually distinct. Direct calculations$^{\text{(}}$%
\footnote{%
This lemma which goes back to Darboux is proved in \cite[lemma 2.4]%
{DoP-HeG1997}.}$^{\text{)}}$ shows that these functions satisfy the shock
wave equation (\ref{F/ SW}).

Let $\eta :\mathcal{Y}\rightarrow \mathbb{CP}_{2}$ the canonical injection.
Since $\eta ^{\ast }\Omega _{\xi }$ may only have poles in $\mathcal{Y}_{\xi
}\cup \mathcal{Y}_{\infty }$, the Stokes formula gives that near $\xi _{\ast
}$, $G=H+L$ where\vspace{-2pt}%
\begin{equation*}
H\left( \xi \right) =\dsum\limits_{z\in \mathcal{Y}_{\xi }}\limfunc{Res}%
\left( \eta ^{\ast }\Omega _{\xi },z\right) \hspace{0.5cm},\hspace{0.5cm}%
L\left( \xi \right) =\dsum\limits_{z\in \mathcal{Y}_{\infty }}\limfunc{Res}%
\left( \eta ^{\ast }\Omega _{\xi },z\right) .
\end{equation*}%
By construction, $\eta ^{\ast }\Omega _{\xi }$ has residue $h_{j}\left( \xi
\right) $ at $z=H_{j}\left( \xi \right) \in \mathcal{Y}_{\xi }$ and it
remains only to know that $L$ is affine in $\xi _{0}$ to prove theorem~\ref%
{T/ caract - a}. The second part of the lemma below is needed in the proof
of theorem~\ref{T/ caract - a - Eff}.

\begin{lemma}
\label{L/ diffaffine}If $W_{\xi _{\ast }}$ is small enough, $L=\Sigma
h_{j}-G $ is affine in $\xi _{0}$. In addition, there is an integer $q$ such
that $L$ is the limit in $\mathcal{O}\left( W_{\xi ^{\ast }}\right) $ of a
continuous one parameter family of $\xi _{0}$-affine functions which are sum
of $q$ mutually distinct shock wave functions.
\end{lemma}

\begin{proof}
With no loss of generality, we assume $\xi _{\ast }=0$ for the proof. For
small complex parameters $\varepsilon $, we consider the homogeneous
coordinates $w^{\varepsilon }=\left( w_{0}+\varepsilon
w_{1}:w_{1}:w_{2}+\varepsilon w_{1}\right) $. For $\varepsilon $ in a
sufficiently small neighborhood of $0$ the intersection of $\mathcal{Y}$
with the zero set of $\Lambda _{\xi }:w\mapsto \xi _{0}+\xi
_{1}w_{1}^{\varepsilon }+w_{2}^{\varepsilon }$ is still generic in the sense
that it is transverse and lies in $\left\{ w_{2}^{\varepsilon }\neq
0\right\} \cap \mathcal{Y}_{\func{reg}}$. Hence, setting $\Omega _{\xi
}^{\varepsilon }=\frac{w_{1}}{w_{0}^{\varepsilon }}\frac{d\Lambda _{\xi
}^{\varepsilon }\left( w\right) }{\Lambda _{\xi }^{\varepsilon }\left(
w\right) }-\frac{w_{1}}{\left( w_{0}^{\varepsilon }\right) ^{2}}%
dw_{0}^{\varepsilon }$, the function $G^{\varepsilon }:\xi \mapsto \frac{1}{%
2\pi i}\int_{\delta }\Omega _{\xi }^{\varepsilon }$ is, on $W_{\xi _{\ast }}$
the sum of $p$ mutually distinct shock wave functions $h_{1}^{\varepsilon
},...,h_{p}^{\varepsilon }$. For generic $\varepsilon $, $\mathcal{Y}$ meets
transversely $\Delta _{\infty }^{\varepsilon }=\left\{ w_{0}^{\varepsilon
}=0\right\} $ and $\mathcal{Y}_{\infty }^{\varepsilon }=\mathcal{Y}\cap
\left\{ w_{0}=0\right\} $ lies in $\mathcal{Y}_{\func{reg}}\cap \left\{
w_{2}^{\varepsilon }\neq 0\right\} $. Hence, \cite[Lemme 2.3.1]{DoP-HeG1997}
implies that $L^{\varepsilon }=\Sigma h_{j}^{\varepsilon }-G^{\varepsilon }$
is affine in $\xi _{0}$. The dependence of $G^{\varepsilon }$ is clearly
holomorphic in $\varepsilon $. The same holds for each $h_{j}^{\varepsilon }$
since what precedes have shown that $h_{j}^{\varepsilon }\left( \xi \right) =%
\limfunc{Res}\left( \eta ^{\ast }\Omega _{\xi }^{\varepsilon },H_{j}\left(
\xi \right) \right) =\frac{1}{2\pi i}\int_{\mathcal{Y}\cap \partial
U_{j}}\eta ^{\ast }\Omega _{\xi }^{\varepsilon }$ where $U_{j}$ is any
sufficiently small neighborhood of $H_{j}\left( \xi \right) $ in $\mathbb{CP}%
_{2}$ whose boundary is smooth and transverse to $\mathcal{Y}$. Hence $%
L^{\varepsilon }$ is holomorphic in $\varepsilon $ and has to be affine in $%
\xi _{0}$ when $\varepsilon =0$.

Let $q$ be the number of points in $\mathcal{Y}_{\infty }$ counted with
their multiplicities~; when $\xi $ is generic, $q$ is either defined by%
\begin{equation}
p-q=\frac{1}{2\pi i}\int_{\gamma }\frac{d\left( \xi _{0}+\xi
_{1}f_{1}+f_{2}\right) }{\xi _{0}+\xi _{1}f_{1}+f_{2}}.  \label{F/ p-q}
\end{equation}%
For sufficiently small generic $\varepsilon $, \cite[Lemme 2.3.1]%
{DoP-HeG1997} gives more precisely that $L^{\varepsilon }=\underset{%
1\leqslant j\leqslant q}{\Sigma }h_{j}^{\varepsilon ,\infty }$ with%
\begin{equation*}
h_{j}^{\varepsilon ,\infty }=-\limfunc{Res}\left( \eta ^{\ast }\Omega _{\xi
},z_{j}^{\varepsilon }\right) =\frac{-\xi _{0}\psi _{j}^{\varepsilon }\left(
0\right) +\psi _{j}^{\varepsilon }{}^{\prime }\left( 0\right) }{1+\xi
_{1}\psi _{j}^{\varepsilon }\left( 0\right) }
\end{equation*}%
where $\mathcal{Y}_{\infty }^{\varepsilon }=\left\{ z_{1}^{\varepsilon
},...,z_{q}^{\varepsilon }\right\} $ and $\psi _{j}^{\varepsilon }\in 
\mathcal{O}\left( U^{\varepsilon }\right) $, $U^{\varepsilon }$ open
neighborhood of $0$ in $\mathbb{C}$, enable to give in the affine
coordinates $\zeta ^{\varepsilon }=(w_{j}^{\varepsilon }/w_{2}^{\varepsilon
})_{j=0,1}$ the set $\mathcal{Y}$ as a graph above $U^{\varepsilon }$~: $%
\mathcal{Y}\cap V_{j}^{\varepsilon }=\left\{ (\zeta _{0}^{\varepsilon }:\psi
_{j}^{\varepsilon }\left( \zeta _{0}^{\varepsilon }\right) :1)~;~\zeta
_{0}^{\varepsilon }\in U^{\varepsilon }\right\} $. Each $h_{j}^{\varepsilon
,\infty }$ is clearly a shock wave function, that is a solution to $h_{\xi
_{1}}=h_{\xi _{0}}h$.
\end{proof}

\noindent \textbf{Remark\textit{.}\ }When $\varepsilon $ goes to a non
generic value, the fact that $L$ is a sum of $q$ shock wave functions may
not be preserved as section~\ref{S/ EffCha} shows.

\subsection{Proof of B.theorem~\protect\ref{T/ caract - a}\label{S/ 3aB}}

Assume that $\gamma $ satisfies (\ref{F/ caract G}) in a connected
neighborhood $W_{\xi _{\ast }}$ of one point $\left( \xi _{0\ast }:\xi
_{1\ast }:1\right) $ of $\mathbb{CP}_{2}$. If $\left( \partial
^{2}G/\partial \xi _{0}^{2}\right) _{\left\vert W_{\xi _{\ast }}\right. }=0$%
, then $\gamma $ satisfies the classical Wermer-Harvey-Lawson moment
condition in $\mathbb{C}_{\xi _{\ast }}^{2}=\mathbb{CP}_{2}\backslash
\left\{ \xi _{0\ast }w_{0}+\xi _{1\ast }w_{1}+w_{2}=0\right\} $ (see \cite[%
cor. 1.6.2]{DoP-HeG1997}) and \cite{WeJ1958c}\cite{HaR-LaB1975} implies that
if $\delta $ is suitably oriented, the polynomial hull of $\delta $ in $%
\mathbb{C}_{\xi _{\ast }}^{2}$ is the unique complex curve $\mathcal{Y}$ of
finite mass of $\mathbb{C}_{\xi _{\ast }}^{2}\backslash \delta $ such that $d%
\left[ \mathcal{Y}\right] =\left[ \delta \right] $.

Assume now $\left( \partial ^{2}G/\partial \xi _{0}^{2}\right) _{\left\vert
W_{\xi _{\ast }}\right. }\neq 0$. Then we can choose a minimal $\mathcal{H}%
=\left\{ h_{1},...,h_{p}\right\} $ in the sense that no proper subset of $%
\mathcal{H}$ satisfy (\ref{F/ caract G}). Although it is not explicitly
mentioned by their authors, the heart of the arguments of \cite[th. II p. 390%
]{DoP-HeG1997} is that $\pm \left[ \delta \right] =d\left[ \mathcal{Y}\right]
$ where $\mathcal{Y}$ is the analytic extension $\mathcal{Y}$ in $\mathbb{CP}%
_{2}\backslash \delta $ of the union $\Gamma $ of the graphs $\Gamma _{j}$
of the functions 
\begin{equation*}
H_{j}:\xi \mapsto \left( 1:h_{j}\left( \xi \right) :-\xi _{0}-\xi
_{1}h_{j}\left( \xi \right) \right) ,1\leqslant j\leqslant p.
\end{equation*}%
This fact, not totally explicit in \cite[p. 264]{HeG1995}, can be recovered 
\textit{a posteriori} by a kind of trick which has been used Poly in \cite%
{DoP-PoJ1995} and is developed later in the proof of theorem~\ref{T/ caract
- c}~: for the curve $\widetilde{\gamma }$ which is the union of $\gamma $
with the boundaries of $\Gamma _{j}$ negatively oriented, one goes back to
the $\mathbb{C}_{\xi _{\ast }}^{2}$-case where $\left( \partial
^{2}G/\partial \xi _{0}^{2}\right) _{\left\vert W_{\xi _{\ast }}\right. }=0$%
. If $d\left[ \mathcal{Y}\right] $ is $-\left[ \delta \right] $ and not $%
\left[ \delta \right] $, then the same arguments which have proved theorem~%
\ref{T/ caract - a}.A would give that the functions $h_{j}$, geometrically
defined as the first coordinates of points of intersection of $\mathcal{Y}$
with generic lines $\Lambda _{\xi }$, should satisfy not only the shock wave
equation $h_{\xi _{1}}=h_{\xi _{0}}h$ but also the "negative" shock wave
equation $h_{\xi _{1}}=-h_{\xi _{0}}h$. As this is impossible, $d\left[ 
\mathcal{Y}\right] =\left[ \delta \right] $.

In both cases, we have found (up to a change of orientation if $\partial
^{2}G/\partial \xi _{0}^{2}$ vanish on $W_{\xi _{\ast }}$) a complex curve $%
\mathcal{Y}$ of finite mass of $\mathbb{CP}_{2}\backslash \delta $ such that 
$d\left[ \mathcal{Y}\right] =\left[ \delta \right] $.

As $\delta $ is smooth, we know from~\cite{HaR-LaB1975} that there is in $%
\delta $ a compact set $\tau $ such that $h^{1}\left( \tau \right) =0$ and
for which each point of $y\in \delta \backslash \tau $ has a neighborhood $%
U_{y}$ where $\overline{\mathcal{Y}}\cap U_{y}$ is a closed bordered
submanifold of $U_{y}$ with boundary $\delta \cap U_{y}$. Lemma~\ref{L/ reg
bord} in section~\ref{S/ UsousE} describes how $\mathcal{Y}$ is near points
of $\tau $. Using the notations and definitions introduced after its proof,
we let $\tau _{2}$ (resp. $\tau ^{\prime }$) be the set of $y$ in $\delta $
where $\overline{\mathcal{Y}}$ has a strong singularity (resp. $%
m_{y}\geqslant 2$) and define $\widetilde{\mathcal{Y}}$ as the abstract
complex curve $\mathcal{Y}\cup \tau ^{\prime }$.

Consider a normalization $\pi :\mathcal{X}\rightarrow \widetilde{\mathcal{Y}}
$~; lemma~\ref{L/ reg bord} implies that $\pi $ is an open mapping. Let $%
\mathcal{Z}$ be the disjoint and abstract union $\mathcal{X}\cup \gamma $.
If $x\in \gamma $, we define a neighborhood of $x$ in $\mathcal{Z}$ as a
subset of $\mathcal{Z}$ which contains a set of the kind $\pi ^{-1}\left(
C_{y,1}^{U}\right) $ where $y=f\left( x\right) $ and $U$ is a neighborhood
of $y$ in $\mathbb{CP}_{2}$. Then $\left( \mathcal{Z},\gamma \right) $ is a
compact metrizable topological manifold with boundary which has finite
2-dimensionnal Hausdorff measure and smooth boundary outside $\sigma
=f^{-1}\left( \tau _{2}\right) $. Since $f$ is an embedding, $h^{1}\left(
\sigma \right) =0$ and $\left( \mathcal{Z},\gamma \right) $ is a manifold
with almost smooth boundary. Moreover, it follows by construction that the
meromorphic extension $F:\mathcal{Z}\rightarrow \overline{\mathcal{Y}}$ of $%
f $ to $\mathcal{X}$ defined by $F\left\vert _{\mathcal{X}}\right. =\pi $ is
a normalization of $\overline{\mathcal{Y}}$ in the sense of lemma~\ref{L/
normalisation2}.$_{\square }$\medskip

\noindent \textbf{Remark. }When $\gamma $ is real analytic, \cite[th. II]%
{HaR-LaB1975} implies that $\mathcal{C}_{y,1}^{U}$ is a manifold with
boundary in the classical sense. So, in that case, $\left( \mathcal{Z}%
,\gamma \right) $ is a classical manifold with boundary.\smallskip

The following proposition which clarifies some results of \cite{DoP-HeG1997}
justifies the fourth remark after theorem~\ref{T/ caract - a}.

\begin{proposition}
\label{L/ caract sol}Assume $\delta $ is connected. Then if $\left( \partial
^{2}G/\partial \xi _{0}^{2}\right) _{\left\vert W_{\xi _{\ast }}\right. }=0$%
, the polynomial hull of $\gamma $ in $\mathbb{C}^{2}=\mathbb{CP}%
_{2}\backslash \left\{ \xi _{0\ast }w_{0}+\xi _{1\ast }w_{1}+w_{2}=0\right\} 
$ has boundary $\pm \left[ \gamma \right] $. If $\left( \partial
^{2}G/\partial \xi _{0}^{2}\right) _{\left\vert W_{\xi _{\ast }}\right.
}\neq 0$ and no proper subset of $\mathcal{H}=\left\{
h_{1},...,h_{p}\right\} $ satisfies (\ref{F/ caract G}), then the analytic
extension $\mathcal{Y}$ in $\mathbb{CP}_{2}\backslash \delta $ of the union $%
\Gamma $ of the graphs $\Gamma _{j}$ of the functions $H_{j}:\xi \mapsto
\left( 1:h_{j}\left( \xi \right) :-\xi _{0}-\xi _{1}h_{j}\left( \xi \right)
\right) $, $1\leqslant j\leqslant p$, is the complex curve which has minimal
volume among complex curves $\mathcal{Z}$ such that $d\left[ \mathcal{Z}%
\right] =\left[ \delta \right] $.
\end{proposition}

\begin{proof}
The preceding proof contains the above conclusion except the minimality of
volume of $\mathcal{Y}$ when $\left( \partial ^{2}G/\partial \xi
_{0}^{2}\right) _{\left\vert W_{\xi _{\ast }}\right. }\neq 0$ and $\mathcal{H%
}$ is minimal. In that case, let $\mathcal{Y}^{\prime }$ be a complex curve
of $\mathbb{CP}_{2}\backslash \delta $ with boundary $\left[ \delta \right] $
and minimal volume. Then $\mathcal{Y}=\mathcal{Y}^{\prime }\cup \mathcal{Z}$
where $\mathcal{Z}$ is a union of compact complex curves of $\mathbb{CP}_{2}$%
. But the intersection of any compact curve with a line $\Lambda _{\xi }$ is
described, for $\xi $ in a neighborhood of a generic $\xi _{\ast }$, as a
finite union of graphs of function $\left( 1:g_{j}:g_{j,2}\right) $ whose
second homogeneous coordinate satisfy the shock wave equation and such that $%
\Sigma g_{j}$ is affine in $\xi _{0}$ (see \cite[section 2]{HeG1995}). Since 
$\mathcal{H}$ is minimal, no such $g_{j}$ belongs to $\mathcal{H}$ and it
appears that $\Gamma $ has to be contained in $\mathcal{Y}^{\prime }$ Hence, 
$\mathcal{Y}\subset \mathcal{Y}^{\prime }$ and, finally, $\mathcal{Y}=%
\mathcal{Y}^{\prime }$.\smallskip
\end{proof}

\subsection{Proof of C.theorem~\protect\ref{T/ caract - a}}

Let $\mathcal{X}$, $\mathcal{Z}$, $\mathcal{D}$ and $g$ be as in the
statement. Let $u\in C^{1}\left( \gamma \right) $ and for $z\in \mathcal{D}%
\backslash \gamma $ set%
\begin{equation*}
\Omega _{z}=u\partial _{\zeta }g_{z}+g_{z}\overline{\theta u}.
\end{equation*}%
If $\mathcal{Y}$ is a smooth domain in $\mathcal{Z}$, the Stokes formula (%
\ref{F/ Stokes}) implies that value of $\mathbf{1}_{\mathcal{Y}}\left(
z\right) \widehat{u}\left( z\right) +\int_{\mathcal{Y}}\overline{\partial }%
\widehat{u}\wedge \partial _{\zeta }g_{z}-\int_{\left( \partial \mathcal{Y}%
\right) \backslash \gamma }u\partial _{\zeta }g_{z}$ does not depend of the
Lipschitz extension $\widehat{u}$ of $u$ to $\mathcal{Z}$ and equals $%
\int_{\gamma \cap \mathcal{Y}}u\partial _{\zeta }g_{z}$ when $\mathcal{X}_{%
\func{sing}}=\varnothing $. Hence we can take it as a definition for $%
\int_{\gamma \cap \mathcal{Y}}u\partial _{\zeta }g_{z}$ in the general case.
Let then $F$ be the function defined for $z\in \mathcal{D}\backslash \gamma $
by 
\begin{equation*}
F\left( z\right) =\frac{2}{i}\int_{\zeta \in \gamma }\Omega _{z}\left( \zeta
\right) ~,~~\Omega _{z}=u\partial _{\zeta }g_{z}+g_{z}\overline{\theta u},
\end{equation*}%
where $g_{z}=g\left( .,z\right) $.

Since $\widetilde{u}$, $u$ and $\theta u$ are continuous, the conclusion of
part~C follows from the lemma below which gives $\widetilde{u}=F\left\vert _{%
\mathcal{X}}\right. $ and $\theta u=\partial \widetilde{u}$ on $\gamma
\backslash \sigma $ if $F\left\vert _{\mathcal{D}\backslash \overline{%
\mathcal{X}}}\right. =0$.

\begin{lemma}
\label{L/ saut pour u et theta}$F_{+}=F\left\vert _{\mathcal{X}}\right. $
and $F_{-}=F\left\vert _{\mathcal{D}\backslash \overline{\mathcal{X}}%
}\right. $ are real valued harmonic functions such that%
\begin{equation}
u=F_{+}-F_{-}~~~~~\&~~~~~\overline{\,\theta u}=\overline{\partial }F_{+}-%
\overline{\partial }F_{-}~~~on~~~\gamma \backslash \sigma .
\label{F/ saut pour u et theta_u}
\end{equation}
\end{lemma}

\begin{proof}
The harmonicity of $F$ is a simple consequence of the properties of $g$. Fix
now $p$ in $\gamma \backslash \sigma $ and in neighborhood $\mathcal{U}$ of $%
p$ in $\mathcal{D}$, a holomorphic chart $\mathcal{U}\rightarrow U$ centered
at $p$~; a hat "$~\widehat{}\,~$" denotes thereafter the coordinate
expression of a function, a form or a set. For $z\in U\backslash \gamma $
let us write $F\left( z\right) =\varphi \left( z\right) +R_{1}\left(
z\right) $ with $\varphi \left( z\right) =-2i\int_{\gamma \cap \mathcal{U}%
}\Omega _{z}$ and $R_{1}$ is smooth on $\mathcal{U}$. If $y$ and $x$ are the
coordinates of $\zeta \in \gamma \cap \mathcal{U}$ and $z\in U\backslash
\gamma $, $\widehat{g}\left( y,x\right) =g\left( \zeta ,z\right) $ can be
written in the form $\widehat{g}\left( y,x\right) =\frac{1}{2\pi }\ln
\left\vert y-x\right\vert +h\left( y,x\right) $ where $h$ is a smooth
function on $U\times U$, harmonic in each its variable. Hence%
\begin{equation*}
\widehat{\varphi }\left( x\right) =\frac{1}{2\pi i}\int_{\widehat{\gamma }%
\cap U}\frac{\widehat{u}\left( y\right) }{y-x}dy+\int_{\widehat{\gamma }\cap
U}\frac{\ln \left\vert y-x\right\vert }{\pi i}\overline{\widehat{\theta u}%
\left( y\right) }+R_{2}\left( x\right)
\end{equation*}%
where $R_{2}\in C^{0}\left( U\right) $. The second integral has no jump
across $\widehat{\gamma }$ and from the classical Sohotsky-Plemelj formula,
we know that the first integral has jump $\widehat{u}$ across $\widehat{%
\gamma }$ in the sense of distribution and pointwise near each regular
boundary point. Since%
\begin{equation*}
\overline{\partial }\widehat{\varphi }\left( x\right) =\frac{1}{2\pi i}d%
\overline{x}\int_{\widehat{\gamma }\cap U}\frac{1}{\overline{x}-\overline{y}}%
\overline{\widehat{\theta u}\left( y\right) }+\overline{\partial }%
R_{2}\left( x\right) ,
\end{equation*}%
the jump of $\overline{\partial }\widehat{\varphi }$ through $\widehat{%
\gamma }$ is likewise$\overline{\text{ }\widehat{\theta u}}$.

In order to check that $F\left( z\right) \in \mathbb{R}$ when $z\in \mathcal{%
D}\backslash \gamma $, we let $L_{\mathcal{X}}$ be the DN-operator of $%
\mathcal{X}$ and we note that since $\partial g_{z}=\left( L_{\mathcal{X}%
}g_{z}\right) \left( \nu ^{\ast }+i\tau ^{\ast }\right) $ and $\theta
u=\left( Lu\right) \left( \nu ^{\ast }+i\tau ^{\ast }\right) $, $-\limfunc{Im%
}F\left( z\right) =\int_{\gamma }\left( u\,\tau g_{z}+g_{z}\tau u\right)
\tau ^{\ast }=\int_{\gamma }d\left( ug_{z}\left\vert _{\gamma }\right.
\right) =0$.
\end{proof}

\subsection{Proof of A.theorem~\protect\ref{T/ caract - c}}

We assume that $\gamma $ is in the\ sense of currents the boundary of a
complex curve $\mathcal{X}$ of $\mathbb{CP}_{2}\backslash \gamma $ for which 
$\left( \gamma ,u,\theta u\right) $ is a restricted DN-datum. Then as in the
proof of theorem~\ref{T/ caract - a}.A, for almost all $\left( \xi _{\ast
0}:\xi _{\ast 1}:1\right) $ in $\mathbb{CP}_{2}$ there exists a neighborhood 
$W_{\xi _{\ast }}$ of $\xi _{\ast }=\left( \xi _{\ast 0},\xi _{\ast
1}\right) $ such that for every $\xi \in W_{\xi _{\ast }}$, $\mathcal{X}\cap
\Delta _{\xi }$ lies in $\mathbb{C}^{2}$ and equals $\Gamma \cap \Delta
_{\xi }$ where $\Gamma $ is the union of the graphs $\Gamma _{j}$ of $%
H_{j}=\left( 1:h_{j}:h_{j,2}\right) :W_{\xi _{\ast }}\rightarrow \mathbb{CP}%
_{2}$, $1\leqslant j\leqslant p$ where $H_{j}$ is holomorphic in $W_{\xi
_{\ast }}$. Since we are concerned only by generic $\xi _{\ast }$, we can
suppose that $\Gamma _{j}=\left\{ \left( \varphi _{j}\left( z_{2}\right)
,z_{2}\right) ~;~z_{2}\in U_{j}\right\} $ where $U_{j}$ is a neighborhood of 
$z_{j,2}^{\ast }=h_{j,2}\left( \xi _{\ast }\right) $ and $\varphi _{j}\in 
\mathcal{O}\left\{ U_{j}\right\} $. The decomposition sought for $\widetilde{%
G}$ in (a1) can be then found in~\cite{HeG1995}. However, this residues
calculus is needed in part~B and we include it here

When $\xi \in W_{\xi _{\ast }}$, $\widetilde{G}_{\ell }\left( \xi _{0}:\xi
_{1}:1\right) $ is the sum of the residues of 
\begin{equation*}
\Lambda _{\ell }=\frac{z_{\ell }}{\xi _{0}+\xi _{1}z_{1}+z_{2}}\Theta _{0}
\end{equation*}%
(for convenience $z_{\ell }=1$ if $\ell =0$) in $\mathcal{X}$. Set $\Theta
_{0}=A_{j}\left( z_{2}\right) dz_{2}$ in each $\Gamma _{j}$ and let us
abbreviate $H_{j}\left( \xi \right) $ in $z_{j}$. Then $z_{j,2}$ is the only
pole of $\Lambda _{\ell }$ in $\Gamma _{j}$. It is a simple one and the
residue of $\Lambda _{\ell }$ at it is 
\begin{equation}
g_{j,\ell }=\frac{z_{j,\ell }~A_{j}\left( z_{j,2}\right) }{\xi _{1}\varphi
_{j}^{\prime }\left( z_{j,2}\right) +1}=z_{j,\ell }~g_{j,0}
\label{F/ 1ere forme des g_j}
\end{equation}%
where $z_{j,\ell }=1$ if $\ell =0$. As $\Gamma _{j}$ is also parametrized by 
$H_{j}$, we can set 
\begin{equation*}
g_{j}=\underset{0\leqslant \ell \leqslant 2}{\Sigma }g_{j,\ell }\,d\eta
_{\ell }~~~\text{on}~~~\Gamma _{j}
\end{equation*}%
and get that $g_{j,1}/g_{j,0}=z_{j,1}=h_{j}$ satisfy (\ref{F/ SW})~; $%
g_{j,2}/g_{j,0}=z_{j,2}=h_{j,2}$ satisfy then $h_{j}\frac{\partial h_{j,2}}{%
\partial \xi _{0}}=\frac{\partial h_{j,2}}{\partial \xi _{1}}$ because $%
h_{j,2}=-\xi _{0}-\xi _{1}h_{j}$. Note that the dependence in $\xi $ of $g$
can be made more clear if the identity $h_{j}-\varphi _{j}\left( -\xi
_{0}-\xi _{1}h_{j}\right) =0$ is used. Indeed, this relation implies%
\begin{equation*}
\left( 1+\xi _{1}\varphi _{j}^{\prime }\right) \partial _{\xi
_{0}}h_{j}+\varphi _{j}^{\prime }=0~~~\&~~~\left( 1+\xi _{1}\varphi
_{j}^{\prime }\right) \partial _{\xi _{1}}h_{j}+h_{j}\varphi _{j}^{\prime }=0
\end{equation*}%
Hence $\varphi _{j}^{\prime }=-\left( \partial _{\xi _{0}}h_{j}\right)
/\left( 1+\xi _{1}\partial _{\xi _{0}}h_{j}\right) $ and 
\begin{equation*}
\frac{1}{1+\xi _{1}\varphi _{j}^{\prime }}=\frac{\partial _{\xi _{1}}h_{j}}{%
\partial _{\xi _{0}}h_{j}}\frac{1+\xi _{1}\partial _{\xi _{0}}h_{j}}{h_{j}}%
=1+\xi _{1}\partial _{\xi _{0}}h_{j}=\partial _{\xi _{0}}h_{j,2}
\end{equation*}%
since $h_{j}$ satisfies $h_{j}\partial _{\xi _{0}}h_{j}=\partial _{\xi
_{1}}h_{j}$. So we have now instead of (\ref{F/ 1ere forme des g_j})%
\begin{equation}
g_{j,\ell }=A_{j}\left( H_{j}\right) \,h_{j,\ell }\,\frac{\partial h_{j,2}}{%
\partial \xi _{0}},~1\leqslant j\leqslant p,~0\leqslant \ell \leqslant 2.
\label{F/ 2eme formule pour les g_j}
\end{equation}%
where $h_{j,0}=1$ for convenience.\smallskip

From the definition we get that when expressed in the affine coordinates $%
\xi $, $g_{j}$ is given by the integral formula%
\begin{equation*}
2\pi i~g_{j}=\left( \int_{\partial \Gamma _{j}}\frac{\Theta _{0}}{\xi
_{0}+\xi _{1}z_{1}+z_{2}}\right) d\xi _{0}+\left( \int_{\partial \Gamma _{j}}%
\frac{z_{1}\Theta _{0}}{\xi _{0}+\xi _{1}z_{1}+z_{2}}\right) d\xi _{1}
\end{equation*}%
from which it is clear that $g_{j}$ is closed.\smallskip

To achieve the proof of part~A, it is enough to remark that (a2) is a direct
consequence of the fact $\func{Re}\Theta _{\ell }=d\widetilde{u_{\ell }}$ is
exact.

\subsection{Proof of B.theorem~\protect\ref{T/ caract - c}}

Assume that the hypothesis of (b1) is true and $\gamma $ is connected.
\medskip

\noindent \textit{Case} $\widetilde{G}\left\vert _{W_{\eta _{\ast }}}\right.
\neq 0$. This mean $\mathcal{H}\neq \varnothing $ when $\mathcal{H}$ is
minimal in the sense that no proper subset of $\mathcal{H}$ gives a
decomposition of $\widetilde{G}$ with the same properties. Let $\Gamma $ be
the union of the graphs $\Gamma _{j}$ of the functions $H_{j}=\left(
1:h_{j}:h_{j,2}\right) $, $1\leqslant j\leqslant p$, where $h_{j,2}=-\xi
_{0}-\xi _{1}h_{j}$. If needed, we can choose another $\xi _{\ast }$ in
order that $\Gamma $ does not meet $\Lambda _{\xi _{\ast }}$ in $\left\{
w_{0}=0\right\} $. Then for any $\xi $ in a neighborhood $\Omega $ of $\xi
_{\ast }$, the $H_{j}\left( \xi \right) $ are mutually distinct and are the
points of $\Gamma \cap L_{\xi }$. Finally, we assume, which it is not a
restriction, that $\Gamma $ has a smooth oriented boundary $\partial \Gamma $%
.

Let $\widetilde{\gamma }$ be the union of $\partial \Gamma $ with opposite
orientation and $\gamma $ and let $\varphi _{\eta }$ be the linear function $%
z\mapsto \eta _{0}+\eta _{1}z_{1}+\eta _{2}z_{2}$. From the hypothesis we
get directly%
\begin{equation*}
\frac{1}{2\pi i}\int_{\gamma }\varphi _{\eta }^{-1}\theta _{\ell
}=\dsum\limits_{1\leqslant j\leqslant p}g_{j,\ell }=\dsum\limits_{1\leqslant
j\leqslant p}g_{j,\ell }h_{j,\ell }
\end{equation*}%
where $h_{j,\ell }=1$ if $\ell =0$. On the other hand, if we set 
\begin{equation}
\Theta _{0}=(\partial _{\xi _{0}}h_{j,2})^{-1}g_{j,0}\,dh_{j,2},~~\text{on ~}%
\overline{\Gamma _{j}},~1\leqslant j\leqslant p  \label{F/ theta0 sur Gamma}
\end{equation}%
and $z_{0}=1$, the residues calculus made in the proof of part~A implies that%
\begin{equation*}
\int_{\partial \Gamma }z_{\ell }\varphi _{\eta }^{-1}\theta
_{0}=\dsum\limits_{1\leqslant j\leqslant p}g_{j,0}h_{j,\ell }~,~1\leqslant
j\leqslant p.
\end{equation*}%
Hence, $\int_{\widetilde{\gamma }}\varphi _{\eta }\theta _{\ell }=0$. As $%
\widetilde{\gamma }$ is contained in the affine space $E_{\xi _{\ast }}=%
\mathbb{CP}_{2}\backslash \Lambda _{\xi _{\ast }}$, we can apply \cite[cor.
4.2 p. 265]{HeG1995} and \cite[prop. 1]{DiT1998b} and get in $E_{\xi _{\ast
}}\backslash \widetilde{\gamma }$ a complex curve $\widetilde{\mathcal{X}}$
of finite volume where $\theta _{0}$ extends weakly in a weakly holomorphic
form $\Theta _{0}$ satisfying$^{\text{(}}$\footnote{%
In \cite{HeG1995}\cite{DiT1998b}, (\ref{F/ prolongement de Theta}) is in
fact obtained only for $\varphi $ smooth on $\mathcal{X}$ and holomorphic in
a neighborhood of $\gamma $ but (\ref{F/ prolongement de Theta}) follows
from this together with (\ref{F/ 2eme formule pour les g_j}) and the residue
relations (very close in spirit to the relations \ref{F/ equation Theta}) 
\begin{equation*}
\frac{1}{2\pi i}\int_{\gamma }z_{1}^{m}\frac{\theta _{0}}{\xi _{0}+\xi
_{1}z_{1}+z_{2}}=\underset{1\leqslant j\leqslant p}{\Sigma }h_{j}^{m}\left(
\xi \right) g_{j,0}\left( \xi \right)
\end{equation*}%
where $h_{j}\left( \xi \right) $, $h_{j,2}\left( \xi \right) $ and $%
g_{j,0}\left( \xi \right) $ are as above~; indeed these relations enable for
generic $\xi $ the computation of $g_{j,0}\left( \xi \right) $ by a
kramerian system and hence imply the smoothness of $\Theta _{0}$ near points
of $\gamma \backslash \sigma $. However, this precision is used not
essentialy in the sequel.}$^{\text{)}}$%
\begin{equation}
\int_{\mathcal{X}}\left( \,\overline{\partial }\varphi \,\right) \wedge
\Theta _{0}=\int_{\mathcal{X}}d\left( \varphi \Theta _{0}\right)
=\int_{\gamma }\varphi \theta _{0}  \label{F/ prolongement de Theta}
\end{equation}%
holds for any $\varphi $ smooth in a neighborhood of $\overline{\mathcal{X}}$
and analytic near $\sigma $.

By construction, $\mathcal{X}=\widetilde{\mathcal{X}}\cup \overline{\Gamma }$
is a complex curve of $\mathbb{CP}_{2}\backslash \gamma $ where $\theta _{0}$
extends as a weakly holomorphic form $\Theta _{0}$, this extension
coinciding in $\Gamma $ with the form defined by (\ref{F/ theta0 sur Gamma}%
). If $\ell =1,2$, the form $\Theta _{\ell }=z_{\ell }\Theta _{0}$ is a
weakly holomorphic extension of $\theta _{\ell }$ to $\mathcal{X}$.

Let $\varepsilon >0$ and let $\mathcal{W}_{\varepsilon }$ an $\varepsilon $%
-neighborhood of $\Lambda _{\xi _{\ast }}$. As $\mathcal{X}_{\varepsilon }=%
\mathcal{X}\backslash \overline{\mathcal{W}_{\varepsilon }}$ lies in the
affine space $E_{\xi _{\ast }}\backslash \mathcal{W}_{\varepsilon }$, \cite[%
th. 4.7]{HaR-LaB1975} implies that 
\begin{equation}
d\left[ \mathcal{X}_{\varepsilon }\right] =n_{\varepsilon }\left[ \gamma %
\right] +\underset{1\leqslant j\leqslant p}{\Sigma }n_{\varepsilon ,j}\left[
\gamma _{\varepsilon ,j}\right]  \label{F/ X_epsilon}
\end{equation}%
where for $1\leqslant j\leqslant p$, $n_{\varepsilon },n_{\varepsilon ,j}\in 
\mathbb{Z}$ and $\gamma _{\varepsilon ,j}$ is the (smooth) boundary of the
smooth (manifold) $\mathcal{W}_{\varepsilon }\cap \Gamma _{j}$. If $%
0<\varepsilon ^{\prime }<\varepsilon $,%
\begin{eqnarray*}
\underset{1\leqslant j\leqslant p}{\Sigma }d\left[ \Gamma _{j}\cap (\mathcal{%
W}_{\varepsilon }\backslash \mathcal{W}_{\varepsilon ^{\prime }})\right]
&=&-d\left[ \mathcal{X}_{\varepsilon }\right] +d\left[ \mathcal{X}%
_{\varepsilon ^{\prime }}\right] \\
&=&\left( -n_{\varepsilon }+n_{\varepsilon ^{\prime }}\right) \left[ \gamma %
\right] +\underset{1\leqslant j\leqslant p}{\Sigma }\left( -n_{\varepsilon
,j}+n_{\varepsilon ^{\prime },j}\right) \left[ \gamma _{\varepsilon ,j}%
\right]
\end{eqnarray*}%
Hence, $n_{\varepsilon }=n_{\varepsilon ^{\prime }}\overset{def}{=}n$ and as
each $\Gamma _{j}\cap \mathcal{W}_{\varepsilon }\backslash \mathcal{W}%
_{\varepsilon ^{\prime }}$ is a smooth manifold with boundary, $%
n_{\varepsilon ,j}=n_{\varepsilon ^{\prime },j}=-1$. Taking now limits in (%
\ref{F/ X_epsilon}) when $\varepsilon $ goes to zero, we get $d\left[ 
\mathcal{X}\right] =n\left[ \gamma \right] $. We suppress from $\mathcal{X}$
any compact component it may have and still denote the result by $\mathcal{X}
$~; note that $\mathcal{X}$ has now to be connected. Since $\mathcal{X}$ is
a complex curve of $\mathbb{CP}_{2}\backslash \gamma $, \cite{ChE1982}
implies that if $n\neq 0$, there is in $\gamma $ a compact set $\sigma $
such that $h^{1}\left( \sigma \right) =0$ and $\left( \,\overline{\mathcal{X}%
},\pm \gamma \,\right) $ is manifold with boundary near points of $\gamma
\backslash \sigma $~; as $\mathcal{X}$ is connected, this implies $n=\pm 1$.
When $n=0$, the structure theorems of Harvey-Shiffman~\cite{HaR1977} implies
that $\mathcal{Z=}\overline{\mathcal{X}}$ is then a complex compact curve of 
$\mathbb{CP}_{2}$~; since $\gamma $ is smooth, $\gamma $ is locally a Jordan
curve of $\mathcal{Z}_{\func{reg}}$ and the points where $\gamma $ may meet
the finite set $\mathcal{Z}_{\func{sing}}$ are only self-intersection points
of $\mathcal{Z}$.\smallskip

Since $\mathcal{X}\cap \Lambda _{\xi }=\left\{ \left( 1:h_{j}\left( \xi
\right) :h_{j,2}\left( \xi \right) \right) ,\text{ }1\leqslant j\leqslant
p\right\} $, the Stokes's formula gives $G=\Sigma h_{j}$ (see proof of
theorem~\ref{T/ caract - a}.A).

Assume that $\left( \partial ^{2}G/\partial \xi _{0}^{2}\right) _{\left\vert
W_{\xi _{\ast }}\right. }\neq 0$. Then $n\neq 0$ because otherwise for $\xi $
closed to $\xi _{\ast }$, the intersection of $\mathcal{X}$ with the line $%
\Lambda _{\xi }$ should have to be the intersection with $\Lambda _{\xi }$
of a compact Riemann surface, namely $\overline{\mathcal{X}}$, which by a
theorem of Reiss would force $\Sigma h_{j}$ to be affine in $\xi _{0}$ (see 
\cite[ch. 5.2]{GrP-HaJ1978Li} or \cite[section 2]{HeG1995}). Reasoning like
in the proof of \ref{T/ caract - a}.B, we also eliminate the possibility $d%
\left[ \mathcal{X}\right] =-\left[ \gamma \right] $ because it would imply
that if $\mathcal{H}=\left\{ h_{1},...,h_{p}\right\} $ is minimal in the
sense that no proper subset of $\mathcal{H}$ satisfy (a1), each $h_{j}$ also
satisfies $h_{j,y}=-h_{j,x}h_{j}$. Hence $n=1$.

When $\left( \partial ^{2}G/\partial \xi _{0}^{2}\right) _{\left\vert W_{\xi
_{\ast }}\right. }=0$, that is when $\Sigma h_{j}$ is affine in $\xi _{0}$,
and when $\overline{\mathcal{X}}$ is not an algebraic curve where $\gamma $
is a slit, $n$ has to be non null and hence is $\pm 1$. Since $\widetilde{G}%
\left\vert _{W_{\eta _{\ast }}}\right. \neq 0$, the minimal $\mathcal{H}$ in
the above sense is not empty and reasoning likewise, we get $d\left[ 
\mathcal{X}\right] =\left[ \gamma \right] $.\smallskip

\noindent \textit{Case} $\widetilde{G}\left\vert _{W_{\eta _{\ast }}}\right.
=0$. This means that the minimal $\mathcal{H}$ is empty. Then, we can apply 
\cite[th. 4.2 p. 264]{HeG1995} in the affine case (see \cite{DiT1998b} for a
generalization and a detailed proof in this case) to get in $\mathbb{C}%
^{2}\backslash \gamma $ (here $\mathbb{C}^{2}$ is the complement of $\Lambda
_{\xi _{\ast }}=\left\{ \xi _{\ast 0}w_{0}+\xi _{1\ast
}w_{1}+w_{2}=0\right\} $ in $\mathbb{CP}_{2}$) a complex curve $\mathcal{X}$
of finite volume where $\theta _{0}$ extends weakly in a weakly holomorphic
form $\Theta _{0}$ satisfying (\ref{F/ prolongement de Theta}).

Since $\theta _{0}$ yields a non zero measure on $\gamma $ which, because $%
\widetilde{G}\left\vert _{W_{\eta _{\ast }}}\right. =0$, is orthogonal to
all polynomials of $\mathbb{C}^{2}\sim \mathbb{CP}_{2}\backslash \Lambda
_{\xi _{\ast }}$, we can apply Bishop~\cite{BiE1963} and \cite{StoG1966}
(Wermer originates and solves this problem \cite{WeJ1958c} for the real
analytic case) to get that $\mathcal{X}$ is the polynomial hull $\widetilde{%
\gamma }$ of $\gamma $ in $\mathbb{C}^{2}$ and $d\left( \pm \left[ \mathcal{X%
}\right] \right) =\left[ \gamma \right] $. \cite{HaR-LaB1975} (see also \cite%
{ChE1982}) implies then that $\gamma $ contains a compact set $\sigma $ such
that $h^{1}\left( \sigma \right) =0$ and $\left( \,\overline{\mathcal{X}}%
,\gamma \,\right) $ is manifold with boundary near points of $\gamma
\backslash \sigma $.\medskip

To prove (b2), we go back to the assumption that (a1) is true and we assume
in addition that $\Theta _{\ell }$ is holomorphic and satisfies (\ref{F/
sans Green1}) for all $c\in H_{1}\left( \mathcal{X}_{\func{reg}}\right) $.
Then, there exists $U_{\ell }\in C^{\infty }\left( \mathcal{X}_{\func{reg}%
}\right) $ such that $\Theta _{\ell }=dV_{\ell }$ on $\mathcal{X}_{\func{reg}%
}~$; since $\Theta _{\ell }$ is a (1,0)-form, $\partial V_{\ell }=dV_{\ell
}=\Theta _{\ell }$. We know from the preceding point that $d\left[ \mathcal{X%
}\right] =n\left[ \gamma \right] $ where $n\in \left\{ 0,1\right\} $, up to
a change of orientation of $\gamma $ when $\widetilde{G}\left\vert _{W_{\eta
_{\ast }}}\right. =0$. Assume at first that $n=1$. As $d$ is elliptic up to
the boundary, $V_{\ell }$ has to be smooth up to the boundary in the
classical sense near points of $\gamma $ outside $\sigma $ and the preceding
equality (\ref{F/ prolongement de Theta}) yields $dv_{\ell }=du_{\ell }$
where $v_{\ell }=V_{\ell }\left\vert _{\gamma }\right. $. Hence, $z_{0}$ is
a point of $\gamma $ where $\left( \,\overline{\mathcal{X}},\gamma \,\right) 
$ is a manifold with boundary, there is a constant $c$ such that $v_{\ell
}=u_{\ell }+c$ near $z_{0}$ in $\gamma $. Since $h^{1}\left( \sigma \right)
=0$ and $dv_{\ell }=du_{\ell }$ is smooth we have $v_{\ell }\left( z\right)
-v_{\ell }\left( z_{0}\right) =\int_{\gamma _{z_{0},z}}du_{\ell }$ where $%
\gamma _{z_{0},z}$ is the positively oriented path of $\gamma $ starting at $%
z_{0}$ and ending at $z$. Hence, $v_{\ell }\left( z\right) -u_{\ell }\left(
z_{0}\right) -c=u_{\ell }\left( z\right) -u_{\ell }\left( z_{0}\right) $ and 
$v_{\ell }\left( z\right) =u_{\ell }\left( z\right) +c$. This implies that $%
U_{\ell }$ is a weakly harmonic extension of $u_{\ell }$. When $n=0$, $%
\overline{\mathcal{X}}$ is two sided locally near points of $\gamma $ and
each local side has the same boundary regularity as in the case $n=1$.
Hence, we can reason as in this case and get that there is a weakly harmonic
extension of $u_{\ell }$ to $\mathcal{X}$.\smallskip

When we assume also that $\int_{\mathcal{X}_{\func{reg}}}\Theta _{\ell
}\wedge \overline{\Theta _{\ell }}<+\infty $, $\Theta _{\ell }$ is
holomorphic in the sense that its pullback to any normalization $\pi :%
\mathcal{Z}\rightarrow \mathcal{X}$ of $\mathcal{X}$ is holomorphic and not
only meromorphic. The isolated singularities that the pullback $V_{\ell }$
of $U_{\ell }$ may have in $\mathcal{Z}$ are removable because $dV_{\ell }$
is smooth. So $U_{\ell }$ is harmonic on $\mathcal{X}$ and $U_{\ell }$ is
the harmonic extension of $u_{\ell }$ to $\mathcal{X}$. The proof is
complete.

\subsection{Proof of theorem~\protect\ref{T/ caract - b}.}

The map $f$ enable us to embed the abstract IDN-problem of theorem~\ref{T/
caract - b} in the projective but concrete frame of theorem~\ref{T/ caract -
c}. So theorem~\ref{T/ caract - b}.A is a direct consequence of \ref{T/
caract - c}.A. \smallskip

For the converse part~B.b1, we apply \ref{T/ caract - c}.B to $\delta
=f\left( \gamma \right) $ and $\theta _{\ell }=\theta u_{\ell }$, $\ell
=0,1,2$. We get in $\mathbb{CP}_{2}\backslash \delta $ an irreducible
complex curve $\mathcal{Y}$ such that $d\left[ \mathcal{Y}\right] =n\left[
\delta \right] $ where $n\in \left\{ 0,1\right\} $, up to a change of
orientation when $\widetilde{G}\left\vert _{W_{\eta _{\ast }}}\right. =0$~;
in addition, each $\theta _{\ell }$ extends weakly to $\mathcal{Y}$ into a
weakly holomorphic (1,0)-form $\Theta _{\ell }^{\mathcal{Y}}$.

When $n=1$, the boundary regularity of $\mathcal{Y}$ mentioned in the proof
of theorem~\ref{T/ caract - c}.B enable to apply readily the construction
made in the proof of theorem~\ref{T/ caract - a}~: adding to $\mathcal{Y}$ a
subset $\sigma ^{\prime }$ of $\gamma $ of zero one dimensional Hausdorff
measure, we get an abstract complex curve $\widetilde{\mathcal{Y}}$ which
can be normalized in the classical sense into an abstract Riemann surface $%
\mathcal{X}$~; $\gamma $ can then be topologically glued to $\mathcal{X}$ so
that $\left( \,\overline{\mathcal{X}},\gamma \,\right) $ becomes a manifold
with almost smooth boundary where the pullback $F$ to $\mathcal{X}$ of the
meromorphic map $\mathbb{CP}_{2}\ni z\mapsto \left( z_{1},z_{2}\right) $
gives a meromorphic extension of $f$ to $\mathcal{X}$. Since the forms $%
\Theta _{\ell }^{\mathcal{Y}}$ are meromorphic on $\mathcal{Y}$, $\Theta
_{\ell }=F^{\ast }\Theta _{\ell }^{\mathcal{Y}}$ is well defined and
meromorphic outside $\mathcal{X}\backslash F^{-1}(\sigma ^{\prime })$ which
has zero length.

When $n=0$, $\mathcal{Z}=\overline{\mathcal{X}}$ is an algebraic curve and
one can use a standard normalization of $\mathcal{Z}$ to get the same kind
of conclusions.\smallskip

The supplementary hypothesis of part~B.b2 force each $\Theta _{\ell }$ to
have only removable singularities. Reasoning like in the proof of \ref{T/
caract - c}.B.b2, (\ref{F/ sans Green1}) implies that $\Theta _{\ell
}=dV_{\ell }=\partial V_{\ell }$ for some harmonic function $V_{\ell }$
smooth up to regular boundary points of $\overline{\mathcal{X}}$ (when $n=0$%
, $\gamma $ cuts locally $\mathcal{X}$ into two domains and this means that
each restriction of $V_{\ell }$ to these domains is smooth up to $\gamma $)
and that there is a constant $c$ such that $U_{\ell }+c$ agrees with $%
u_{\ell }$ on $\gamma $.$_{\square }$

\section{Characterizations, effective or affine}

\subsection{Explicit integro-differential characterization\label{S/ EffCha}}

In this section where th.~\ref{T/ caract - a - Eff} is proved, $\left( \xi
_{0},\xi _{1}\right) $ is replaced by the simpler $\left( x,y\right) $ and
we reason in a neighborhood of $\left( 0,0\right) $. $\mathcal{O}_{0}$ is
the set of one variable holomorphic functions near $0$ is denoted by, $%
\mathcal{O}_{0}\otimes \mathbb{C}_{d}\left[ Z\right] $ stands for the set of
polynomials of degree at most $d$ with independent variable $Z$ and
coefficients in $\mathcal{O}_{0}$. An element of $\mathcal{O}_{0}\otimes 
\mathbb{C}_{d}\left[ X\right] $ (resp. $\mathcal{O}_{0}\otimes \mathbb{C}_{d}%
\left[ Y\right] $) should be think as a function of the type $\left(
x,y\right) \mapsto \underset{0\leqslant j\leqslant d}{\Sigma }\lambda
_{j}\left( x\right) y^{j}$ (resp. $\left( x,y\right) \mapsto \underset{%
0\leqslant j\leqslant d}{\Sigma }\lambda _{j}\left( y\right) x^{j}$) where
each $\lambda _{j}\in \mathcal{O}_{0}$.

If $h$ is differentiable, the derivative of $h$ with respect to one of its
variable $u$ is denoted $h_{u}$. If $U$ is an open set in $\mathbb{C}^{n}$, $%
\mathcal{O}\left( U\right) $ is the space of holomorphic functions in $U$~;
if $h\in \mathcal{O}\left( U\right) $ and $0\in U$, we set $%
h_{x^{-0},y^{0}}=h$ and in any simply connected neighborhood of $0$ in $U$,
we denote by $h_{x^{-\alpha -1},y^{\beta }}$ ($\alpha ,\beta \in \mathbb{N}$%
) the function which vanishes at $0$ and satisfies $\partial \left(
h_{x^{-\alpha -1},y^{\beta }}\right) /\partial x=h_{x^{-\alpha },y^{\beta }}$%
~; $h_{x,^{\alpha }y^{-\beta -1}}$ is defined similarly.

A two variable function $h$ is called a \textit{shock wave function} on a
domain $D$ of $\mathbb{C}^{2}$ if it is holomorphic and satisfies $%
h_{y}=h_{x}h$ on $D$.

A $p$\textit{-algebro\"{\i}de} function on $D$ is a $p$-uple $h=\left(
h_{1},...,h_{m}\right) $ of functions from $D$ to $\mathbb{C}$ for which one
can find $p$ holomorphic functions $a_{0},,...,a_{p-1}$ in $D$ such that for 
$z\in D$, $h_{1}\left( z\right) ,...,h_{p}\left( z\right) $ are the roots,
with multiplicities, of the polynomial $T^{h}=X^{p}+a_{p-1}\left( z\right)
X^{p-1}+\cdots +a_{1}\left( z\right) X+a_{0}\left( z\right) $.

A $p$\textit{-multivaluate shock wave function} on $D$ is a $p$-algebro\"{\i}%
de function $h$ on $D$ such that $\Sigma =\left\{ \func{Discr}%
T^{h}=0\right\} $ is a hypersurface of $D$ and for any $z_{\ast }\in
D\backslash \Sigma $, the holomorphic functions $h_{1}^{\ast
},..,h_{p}^{\ast }$ which near $z_{\ast }$ describe the roots of $T_{z}^{h}$
are non null shock wave functions~; the first symmetric function of $T^{h}$,
that is the sum of the roots of $T^{h}$, is called the \textit{trace} of $%
T^{h}$ of $h$. Traces of $p$-multivaluate shock wave functions are called $p$%
\textit{-shock wave functions}.

%TCIMACRO{\TeXButton{Reprise}{\renewcommand{\theX}{\arabic{X}}}}%
%BeginExpansion
\renewcommand{\theX}{\arabic{X}}%
%EndExpansion

\begin{lemma}
\label{L/ EqDiffSym}Consider $p$ mutually distinct functions $%
h_{1},...,h_{p} $ holomorphic in a domain $D$ of $\mathbb{C}^{2}$. Then,
each $h_{j}$ is a shock wave function if and only the functions $\sigma
_{k}=\left( -1\right) ^{k}\underset{1\leqslant j_{1}<\cdots <j_{k}\leqslant p%
}{\Sigma }h_{j_{1}}...h_{j_{k}}$ satisfy the following system%
\begin{equation}
\sigma _{p}\sigma _{1,x}+\sigma _{p,y}=0\hspace{0.5cm}\&\hspace{0.5cm}\sigma
_{k}\sigma _{1,x}+\sigma _{k,y}=\sigma _{k+1,x}~,~~1\leqslant k\leqslant p-1.
\label{F/ EqSym1}
\end{equation}
\end{lemma}

\begin{proof}
Set $T=X^{p}+\underset{1\leqslant k\leqslant p}{\Sigma }\sigma _{k}X^{p-k}$
and $T^{\prime }=\partial T/\partial X$. The relations $0=\left( Th\right)
_{x}=\left( Th\right) _{y}$ and $h^{p}=-\underset{1\leqslant k\leqslant p}{%
\Sigma }\sigma _{k}h^{p-k}$ yields $(T^{\prime }h)\left( h_{y}-hh_{x}\right)
=Sh$ with%
\begin{equation*}
S=\underset{1\leqslant k\leqslant p-1}{\Sigma }\left[ \sigma _{k+1,x}-\sigma
_{k}\sigma _{1,x}-\sigma _{k,y}\right] X^{p-k}-\left( \sigma _{p}\sigma
_{1,x}+\sigma _{p,y}\right)
\end{equation*}%
Since $\deg S\leqslant p-1$, the fact that each $h_{j}$ is a shock wave
function implies that the coefficients of $S$ vanish in a non empty open set
and thus in the domain $D$. If $S=0$ at every point of $D$, then each $h_{j}$
verifies $h_{y}-h_{x}h=0$ in the domain $D$ because $T^{\prime
}h_{j}\not\equiv 0$.
\end{proof}

\begin{proposition}
\label{P/ ResEqDiffSym}Let $D$ be a simply connected domain of $\mathbb{C}%
^{2}$ containing $0$, $\Delta $ its image by the projection $\left(
x,y\right) \mapsto y$ and and $H\in \mathcal{O}\left( D\right) $. When $u$
is differentiable, we set%
\begin{equation*}
\mathcal{D}_{H}u=e^{H_{x,y^{-1}}}\partial (ue^{-H_{x,y^{-1}}})/\partial y%
\hspace{0.5cm}\&\hspace{0.5cm}\mathcal{L}_{H}u=\left( \mathcal{D}%
_{H}u\right) _{x^{-1}}
\end{equation*}%
The following two assertions are equivalent

1. $H$ is a $p$-shock wave function in $D$.

2. There exists $\lambda _{1},...,\lambda _{p-1}\in \mathcal{O}\left( \Delta
\right) $ such that for $\widetilde{\lambda }_{j}\left( x,y\right) =\lambda
_{j}\left( y\right) $, $1\leqslant j\leqslant p-1$, 
\begin{gather}
\mathcal{D}_{H}\mathcal{L}_{H}^{p-1}H=\mathcal{D}_{H}\mathcal{L}_{H}^{p-2}%
\widetilde{\lambda }_{1}+\cdots +\mathcal{D}_{H}\mathcal{L}_{H}^{0}%
\widetilde{\lambda }_{p-1}  \label{F/ EqNu1} \\
\func{Discr}T_{z}\not\equiv 0  \label{F/ EqNu2}
\end{gather}%
where $T_{z}=X^{p}+\underset{1\leqslant k\leqslant p}{\Sigma }s_{k}\left(
z\right) X^{p-k}$ with%
\begin{equation}
s_{k}=-\mathcal{L}_{H}^{k-1}H+\mathcal{L}_{H}^{k-2}\widetilde{\lambda }%
_{1}+\cdots +\mathcal{L}_{H}^{0}\widetilde{\lambda }_{k-1},~~1\leqslant
k\leqslant p.  \label{F/ Sym}
\end{equation}

More precisely, in case (2) is true, $T$ determines a $p$-multivaluate shock
wave function with trace $H$. Conversely, if $H$ is the trace of a $p$%
-multivaluate shock wave function $T$, the $p$ holomorphic functions which
near a point $z_{\ast }$ in $\left\{ \func{Discr}T\neq 0\right\} $ describes
the roots of $T_{z}$ have symmetric functions $\left( -1\right) ^{k}s_{k}$,$%
1\leqslant k\leqslant p$ which satisfy (\ref{F/ EqNu1}) and (\ref{F/ Sym}).
\end{proposition}

\begin{proof}
1) Assume $H$ is a $p$-shock wave function in $D$. Then, $H$ is the first
symmetric function of some $T\in \mathcal{O}\left( D\right) \otimes \mathbb{C%
}_{p}\left[ Z\right] $, $\Sigma =\left\{ \func{Discr}T=0\right\} \neq D$ and
for any fixed $z_{\ast }\in D\backslash \Sigma $, $\deg T_{z_{\ast }}=p$ and
the holomorphic functions $h_{1}^{\ast },..,h_{p}^{\ast }$ which near $%
z_{\ast }$ describe the roots of $T_{z}$ are non null shock wave functions
with no common value on a sufficient small convex neighborhood $W=U\times V$
of $z_{\ast }$. For $k\in \left\{ 1,...,p\right\} $ and on $W$, set $\rho
_{k}=\sigma _{k}e^{-H_{x,y^{-1}}}$ where $\sigma _{k}$ is defined in lemma~%
\ref{L/ EqDiffSym}. Then (\ref{F/ EqSym1}) implies that $\rho _{p,y}=0$ and
that for $k\in \left\{ 1,...,p-1\right\} $, 
\begin{equation*}
\left( \rho _{k+1}e^{H_{x,y^{-1}}}\right) _{x}=-H_{x}\sigma _{k}+\sigma
_{k,y}=e^{H_{x,y^{-1}}}\frac{\partial }{\partial y}\sigma
_{k}e^{-H_{x,y^{-1}}}=e^{H_{x,y^{-1}}}\rho _{k,y}
\end{equation*}%
which yields $\lambda _{k}\in \mathcal{O}\left( V\right) $ such that%
\begin{equation*}
\rho _{k+1}\left( x,y\right) e^{H_{x,y^{-1}}}=\left[ e^{H_{x,y^{-1}}}\rho
_{k,y}\right] _{x^{-1}}+\lambda _{k}\left( y\right) .
\end{equation*}%
Since $\sigma _{1}=-H$, we get $e^{-H_{x,y^{-1}}}\left[ e^{H_{x,y^{-1}}}\rho
_{1,y}\right] _{x^{-1}}=-e^{-H_{x,y^{-1}}}\mathcal{L}_{H}H$. Setting $%
\widetilde{\lambda }_{0}=-H$, we obtain $\rho _{2}=e^{-H_{x,y^{-1}}}(%
\mathcal{L}_{H}\lambda _{0}+\widetilde{\lambda }_{1})$ and a straightforward
finite recurrence gives (\ref{F/ Sym}) with $\left( s_{k}\right) =\left(
\sigma _{k}\right) $. In particular, $k=p$ yields (\ref{F/ EqNu1}), because $%
\rho _{p,y}=0$~; the discriminant of $T_{z}$ don't vanish in $W$ because $%
h_{1},...,h_{p}$ have no common value. Since (\ref{F/ Sym}) also reads $%
\sigma _{k+1}=\mathcal{L}_{H}\sigma _{k}+\widetilde{\lambda }_{k}$ we obtain
that $\lambda _{k}=\sigma _{k+1}\left( 0,.\right) $ on $V$. Hence, $\lambda
_{1},...,\lambda _{p}$ do not depend on $z_{\ast }$ so that they are well
defined holomorphic function on $\Delta $.

2) Assume now (2) is true. We only have to check that $T=X^{p}-\underset{%
1\leqslant k\leqslant p}{\Sigma }s_{k}X^{p-k}$ is actually a $p$%
-multivaluate shock wave function. Formulas (\ref{F/ Sym}) also read $%
s_{k+1}=\mathcal{L}_{H}s_{k}+\widetilde{\lambda }_{k}$ for $1\leqslant
k\leqslant p-1$ and (\ref{F/ EqNu1}) means that $\mathcal{D}_{H}s_{p}=0$.
Hence $0=s_{p,y}-H_{x}s_{p}=s_{p,y}-s_{1,x}s_{p}$ and $s_{k+1,x}=\mathcal{D}%
_{H}s_{k}=s_{k,y}-H_{x}s_{k}=s_{k,y}+s_{1,x}s_{k}$. So, if $z_{\ast }\in D$
is outside $\Sigma =\left\{ \func{Discr}T=0\right\} $, lemma~\ref{L/
EqDiffSym} implies that the holomorphic functions $h_{1},...,h_{p}$ which
near $z_{\ast }$ describe the roots of $T_{z}$ are mutually distinct shock
wave function.
\end{proof}

The following describes $p$-shock wave functions which are affine in $x$.

\begin{proposition}
\label{P/ OndeAffine}Let $D$ be a simply connected domain of $\mathbb{C}^{2}$
containing $0$, $a,b\in \mathcal{O}\left( D\right) $ and $H=x\otimes
a+1\otimes b$. Then, $H$ is a $p$-shock wave function, if and only if there
exists $Q_{0},Q_{1}\in \mathbb{C}_{p-1}\left[ Y\right] $ such that 
\begin{gather}
a=\frac{Q_{1}}{1-Q_{1,y^{-1}}}~~~\&~~~b=\frac{Q_{0}}{1-Q_{1,y^{-1}}}
\label{F/ a et b} \\
\func{Discr}\left( 1-Q_{1,y^{-1}}\right) \neq 0,~  \label{F/ Qgen}
\end{gather}%
In addition, when (\ref{F/ a et b}) and (\ref{F/ Qgen}) are satisfied, the
decomposition of $H$ in elementary fractions gives $H$ as a sum of shock
wave functions.
\end{proposition}

\begin{proof}
Assume $H=\underset{1\leqslant j\leqslant p}{\Sigma }h_{j}$ where $h_{1}$,
...,$h_{p}$ are mutually distinct shock wave functions~; with proposition~%
\ref{P/ ResEqDiffSym} notations, $\left( -1\right) ^{1}s_{1},...,\left(
-1\right) ^{p}s_{p}$ are the symmetric functions of $h_{1}$, ...,$h_{p}$ and
satisfy the relations $\mathcal{D}_{H}s_{p}=0$ and $s_{k+1}=\mathcal{L}%
_{H}s_{k}+\widetilde{\lambda }_{k}$, $1\leqslant k\leqslant p-1$. There
exist sequences of holomorphic functions, $\left( \lambda _{j,k}\right)
_{\,k\in \mathbb{N}}$, $\left( a_{k}\right) _{k\in \mathbb{N}}$ and $\left(
b_{k}\right) _{k\in \mathbb{N}}$, each satisfying the recurrence $%
u_{k+1}=u_{k}^{\prime }-au_{k}$, such that 
\begin{equation*}
\mathcal{L}_{H}^{k}\widetilde{\lambda }_{j}=\frac{x^{k}}{k!}\otimes \lambda
_{j,k}~,\hspace{0.5cm}\mathcal{L}_{H}^{k}H=\frac{x^{k+1}}{\left( k+1\right) !%
}\otimes a_{k}+\frac{x^{k}}{k!}\otimes b_{k}~,~~k\in \mathbb{N}
\end{equation*}%
Hence, (\ref{F/ EqNu1}) yields the vanishing of the $x$-polynomial%
\begin{equation*}
\frac{x^{p}}{p!}\otimes a_{p}+\frac{x^{p-1}}{\left( p-1\right) !}\otimes
b_{p}-\dsum\limits_{1\leqslant j\leqslant p-1}\frac{x^{p-1-j}}{\left(
p-1-j\right) !}\otimes \lambda _{j,p-j}
\end{equation*}%
and ensures $a_{p}=b_{p}=\lambda _{j,p-j}$, $1\leqslant j\leqslant p-j$. So,
one can find $Q_{1},Q_{0}\in \mathbb{C}_{p-1}\left[ Y\right] $ and $\Lambda
_{j}\in \mathbb{C}_{p-j-1}\left[ Y\right] $, $1\leqslant j\leqslant p-1$,
such that with $A=a_{y^{-1}}$%
\begin{equation*}
a_{k}=Q_{1}^{\left( k\right) }e^{A}~,~~b_{k}=Q_{0}^{\left( k\right)
}e^{A},~~\lambda _{j,k}=\Lambda _{j}^{\left( k\right) }e^{A},~k\in \mathbb{N}%
.
\end{equation*}%
Thus $a=a_{0}=Q_{1}e^{A}$, $1-e^{-A}=Q_{1,y^{-1}}$ and hence, (\ref{F/ a et
b}).

Conversely, assume that (\ref{F/ a et b}) and (\ref{F/ Qgen}) are true for
some $Q_{0},Q_{1}\in \mathbb{C}_{p-1}\left[ Y\right] $. Then, the
decomposition in elementary fraction of $H$ is $\Sigma h_{j}$ where $%
h_{j}\left( x,y\right) =\frac{q_{j}x+c_{j}}{1-q_{j}y}$, $1\leqslant
j\leqslant p$ and $q_{1}^{-1},...,q_{p}^{-1}$ are the roots of $%
1-Q_{1,y^{-1}}$. It is quite evident that $h_{1},...,h_{p}$ are mutually
distinct shock wave functions.
\end{proof}

\noindent \textbf{Remark 1. }If $x\otimes a+1\otimes b$ is the sum of $p$%
-mutually distinct shock wave functions $h_{1}$, ...,$h_{p}$, each $h_{j}$
is algebraic since for $1\leqslant k\leqslant p$, (\ref{F/ Sym}) yields $%
s_{k}\left( x,y\right) =\left[ 1-Q_{1y^{-1}}\left( y\right) \right]
^{-1}S_{k}\left( x,y\right) $ where $S_{k}$ is the polynomial defined by%
\begin{equation*}
S_{k}=-\frac{x^{k}\otimes Q_{1}^{\left( k-1\right) }}{k!}-\frac{%
x^{k-1}\otimes Q_{0}^{\left( k-1\right) }}{\left( k-1\right) !}%
+\dsum\limits_{1\leqslant j\leqslant k-1}\frac{x^{k-1-j}\otimes \Lambda
_{j}^{\left( k-1-j\right) }}{\left( k-1-j\right) !}.
\end{equation*}

\noindent \textbf{Remark 2. }When $Q_{1}$ is allowed to be non generic, $%
1-Q_{1,y^{-1}}$ can only be written in the form $1-Q_{1,y^{-1}}=\underset{%
1\leqslant j\leqslant m}{\Pi }\left( 1-q_{j}y\right) ^{\alpha _{j}}$ with $%
\alpha _{1},...,\alpha _{p}\in \mathbb{N}^{\ast }$. Hence, if $Q_{0}\in 
\mathbb{C}_{p-1}\left[ Y\right] $ and $\left( a,b\right) $ is defined by (%
\ref{F/ a et b}), there is constants $c_{j,\ell }$ such that 
\begin{equation*}
H\overset{def}{=}x\otimes a+1\otimes b=\dsum\limits_{1\leqslant j\leqslant
m}\dsum\limits_{1\leqslant \ell \leqslant \alpha _{j}}h_{j,\ell }
\end{equation*}%
where $h_{j,\ell }=\frac{q_{j}x}{1-q_{j}y}+\frac{c_{j,\ell }}{\left(
1-q_{j}y\right) ^{\ell }}$, $1\leqslant j\leqslant m$, $1\leqslant \ell
\leqslant \alpha _{j}$. Each $h_{j,\ell }$ is now a \textit{generalized
shock wave function} in the sense it is a solution of the equation%
\begin{equation*}
h_{y}-h_{x}h=\left( \ell -1\right) \kappa h_{x}^{\ell +1}
\end{equation*}%
where $\kappa \in \mathbb{C}$ is equal to $c_{j,\ell }/q_{j}^{\ell }$.

Generalized shock wave functions arise when an affine function $H=x\otimes
a+1\otimes b$ has coefficient $a$ and $b$ given by (\ref{F/ a et b}) not
constraint to (\ref{F/ Qgen}). In that case, $H$ is clearly a limit of $p$%
-shock wave functions. The lemma below, which is an elementary consequence
of prop.~\ref{P/ OndeAffine}, proves the converse and so, shows that
generalized shock wave functions occur naturally.

\begin{lemma}
\label{L/ limiteOCaffine}Let $\left( H^{t}\right) _{t\in T}=\left( x\otimes
a^{t}+1\otimes b^{t}\right) _{t\in T}$ be a continuous family of holomorphic
affine functions in a simply connected domain $D$ such that the set $T_{%
\func{reg}}$ of parameters $t$ for which $H_{t}$ is a $p$-shock wave
function is dense in $T$, then there exists in $\mathbb{C}_{p-1}\left[ Y%
\right] $ continuous family of holomorphic polynomials $\left(
Q_{1}^{t}\right) $ and $\left( Q_{0}^{t}\right) $ such that for any $t\in T$%
, $a^{t}=Q_{1}^{t}\left( 1-Q_{1,y^{-1}}^{t}\right) ^{-1}$ and $%
b^{t}=Q_{1,0}^{t}\left( 1-Q_{1,y^{-1}}^{t}\right) ^{-1}$. Hence, $H^{t}$ is
a $p$-shock (resp. $p$-generalized) shock wave function when $\func{Discr}%
\left( 1-Q_{1,y^{-1}}^{t}\right) \neq 0$ (resp. $\func{Discr}\left(
1-Q_{1,y^{-1}}\right) =0$).
\end{lemma}

\subsubsection{Proof of theorem~\protect\ref{T/ caract - a - Eff}.}

\textbf{1)} Assume that there is an open Riemann surface $\mathcal{X}$ such
that $\overline{\mathcal{X}}=\mathcal{X}\cup \gamma $ is a manifold with
almost smooth boundary. Then, theorem~\ref{T/ caract - a} gives that almost
all point $\left( \xi _{0\ast },\xi _{1\ast },1\right) $\ of $\mathbb{CP}%
_{2} $ has a neighborhood $W_{\xi _{\ast }}$ for which one can find an
integer $p$ and $p$ mutually distinct shock wave functions $h_{1},...,h_{p}$
on $W_{\xi _{\ast }}$ such that $L=\Sigma h_{j}-G$ is affine in $\xi _{0}$.
Set $s_{k}=\left( -1\right) ^{k}\underset{1\leqslant j_{1}<\cdots
<j_{k}\leqslant p}{\Sigma }h_{j_{1}}...h_{j_{k}}$ and $T=X^{p}+\underset{%
1\leqslant k\leqslant p}{\Sigma }s_{k}X^{p-k}$. Then $\func{Discr}%
T\not\equiv 0$ and proposition~\ref{P/ ResEqDiffSym} implies that property
(1) of th.~\ref{T/ caract - a - Eff} holds.

Assume now that $\widetilde{p}$ is the least integer $q$ such that there is
a $q$-multivaluate shock wave function whose trace differs from $G\left\vert
_{W_{\xi _{\ast }}}\right. $ only by a $\xi _{0}$-affine function. Let $%
\widetilde{T}$ be a $\widetilde{p}$-multivaluate shock wave function with
trace $\widetilde{H}$ such that $\widetilde{L}=\widetilde{H}-G\left\vert
_{W_{\xi _{\ast }}}\right. $ is affine in $\xi _{0}$. Let $\widetilde{h}%
_{1},...,\widetilde{h}_{\widetilde{p}}$ be the holomorphic function on $%
W_{\xi _{\ast }}$ which describes the roots of $\widetilde{T}$. Then $\{%
\widetilde{h}_{1},...,\widetilde{h}_{\widetilde{p}}\}$ is minimal in the
sense of the fourth remark below theorem~\ref{T/ caract - a}.

When $\widetilde{p}\geqslant 1$, this remark says that $\mathcal{X}$ is a
normalization of the analytic extension $\mathcal{Y}$ in $\mathbb{CP}%
_{2}\backslash f\left( \gamma \right) $ of the union of the graphs of the
functions $\left( 1:\widetilde{h}_{j}:-\xi _{0}-\xi _{1}\widetilde{h}%
_{j}\right) $, $1\leqslant j\leqslant \widetilde{p}$, and for any $\xi \in
W_{\xi _{\ast }}$, the intersection of $\mathcal{Y}$ with the projective
lines $\xi _{0}w_{0}+\xi _{1}w_{1}+w_{2}=0$ is\vspace{-5pt}%
\begin{equation*}
\left\{ \left( 1:\widetilde{h}_{j}\left( \xi \right) :\xi _{0}-\xi _{1}%
\widetilde{h}_{j}\left( \xi \right) \right) ~;~1\leqslant j\leqslant 
\widetilde{p}\right\} .\vspace{-5pt}
\end{equation*}%
Since by prop~\ref{T/ unicite de Y} $\mathcal{Y}$ is uniquely determined by $%
\left( \gamma ,f\right) $, each $\widetilde{h}_{j}$ and so, each symmetric
function $\widetilde{s}_{k}$ of $\widetilde{h}_{1},...,\widetilde{h}_{k}$,
is uniquely determined by $\left( \gamma ,f\right) $. If $\lambda
_{1},...,\lambda _{\widetilde{p}}$ are any one variable holomorphic function
such that assertion (2) of prop.~\ref{P/ ResEqDiffSym} holds, with $%
\widetilde{H}$ instead of $H$, then $\lambda _{k-1}\left( \xi _{1}\right) =%
\widetilde{s}_{k}\left( 0,\xi _{1}\right) $, $1\leqslant k\leqslant p-1$.
Hence, $\lambda _{1},...,\lambda _{\widetilde{p}}$ are uniquely determined
by $\left( \gamma ,f\right) $.

Lemma~\ref{L/ diffaffine} implies now that $L$ is affine in $\xi _{0}$ and
is, for some integer $q$, the limit of a continuous one parameter family of $%
\xi _{0}$-affine $q$-shock wave functions. Thanks to lemma~\ref{L/
limiteOCaffine}, this implies that $L\left( \xi _{0},\xi _{1}\right) =\xi
_{0}\frac{\alpha ^{\prime }\left( \xi _{1}\right) }{1-\alpha \left( \xi
_{1}\right) }+\frac{\beta \left( \xi _{1}\right) }{1-\alpha \left( \xi
_{1}\right) }$ where $\alpha \in \mathbb{C}_{q}\left[ \xi _{1}\right] $
vanish at zero and $b\in \mathbb{C}_{q-1}\left[ \xi _{1}\right] $.\medskip

\textbf{2)} Assume now that property (1) of th.~\ref{T/ caract - a - Eff}
holds. Set $T=X^{p}+\underset{1\leqslant k\leqslant p}{\Sigma }s_{k}X^{p-k}$
where now $s_{k}$ are defined by (\ref{F/ Sym}). Proposition~\ref{P/
ResEqDiffSym} implies then that $T$ is determines a $p$-multivaluate shock
wave function whose trace is $G+L$. Since $L$ is affine in $\xi _{0}$, $G$
satisfies the hypothesis of theorem~\ref{T/ caract - a}.B. As $G$ is not
affine in $\xi _{0}$ by hypothesis, $\gamma $, with its given orientation,
is the boundary of an open Riemann surface $\mathcal{X}$ with the sought
properties. Hence, (4) has to be satisfied from the direct part of theorem~%
\ref{T/ caract - a - Eff}.

\subsubsection{Particular case of theorem~\protect\ref{T/ caract - a - Eff}.}

For minimal $p$ equal to $2$, it turns out that theorem~\ref{T/ caract - a -
Eff} says that $\gamma $ bounds almost smoothly an open Riemann surface $%
\mathcal{X}$ where $f$ extends meromorphically if and only if for some
constants $c,\alpha _{1},\alpha _{2},\beta _{0},\beta _{1}$ 
\begin{equation}
\mathcal{D}_{G}\mathcal{L}_{G}G=cC+\alpha _{1}A_{1}+\alpha _{2}A_{2}+\beta
_{0}B_{0}+\beta _{1}B_{1}  \label{F/ eq p=2b}
\end{equation}%
where%
\begin{eqnarray*}
C &=&\left( G_{x}-g_{x}\right) \mathrm{e}^{g_{x,y^{-1}}} \\
A_{1} &=&y\mathcal{D}_{G}\mathcal{L}_{G}G+\left( x\mathcal{L}_{G}G\right)
_{x}~~,~~~A_{2}=y^{2}\mathcal{D}_{G}\mathcal{L}_{G}G+\left( 2xy\mathcal{L}%
_{G}G+x^{2}G\right) _{x} \\
B_{0} &=&\mathcal{D}_{G}G-\mathcal{D}_{G}g~~,~~B_{1}=yB_{0}+g-G
\end{eqnarray*}

\subsection{Affine characterization\label{S/ AffCha}}

A characterization for an affine presentation is possible but very special
data have to be selected. Assume $\mathcal{X}$ is a Riemann surface with
almost smooth boundary $\gamma $. When $u\in C^{\infty }\left( \gamma
\right) $, a straightforward computation gives that $d^{c}\widetilde{u}%
=\left( Nu\right) \tau ^{\ast }$ as forms of $\gamma \backslash \overline{%
\mathcal{X}}_{\func{sing}}$. So, lemma~\ref{L/ Stokes} implies that%
\begin{equation}
\int_{\gamma }\left( Nu\right) \,\tau ^{\ast }=\int_{\mathcal{X}}dd^{c}%
\widetilde{u}=0.\vspace{-3pt}  \label{F/ obvious}
\end{equation}%
Hence $\left( Nu\right) \tau ^{\ast }$ has a primitive $v$ on $\gamma $ and
the holomorphic extension to $\mathcal{X}$ of $h=u+iv$ is equivalent to the
moments condition\vspace{-2pt}%
\begin{equation}
\forall \Psi \in H^{1,0}(\overline{\mathcal{X}}),~\int_{\gamma }h\Psi =0.%
\vspace{-3pt}  \label{F/ mom}
\end{equation}%
If $H\in \mathcal{O}(\,\overline{\mathcal{X}}\,)^{N+1}$ is such that $%
h=H\left\vert _{\gamma }\right. $ embeds $\gamma $ in $\mathbb{C}^{N+1}$, $%
\mathcal{Y}=H\left( \mathcal{X}\right) \backslash h\left( \gamma \right) $
is a complex curve of $\mathbb{C}^{N+1}\backslash h\left( \gamma \right) $
which has finite volume and has boundary $h\left( \gamma \right) $ and $%
\gamma $ has to satisfy the Harvey-Lawson-moments condition\vspace{-3pt}%
\begin{equation}
\forall k_{0},...,k_{N}\in \mathbb{N},~~\int_{\gamma
}h_{0}^{k_{0}}...h_{N}^{k_{N}}dh_{0}=0\vspace{-3pt}  \label{F/ mom affine}
\end{equation}%
which by the way contains also (\ref{F/ obvious}).\vspace{-3pt}

\begin{proposition}
\label{T/ affine}Let $u_{0},...,u_{N}\in C^{\infty }\left( \gamma \right) $\
and let $v_{\ell }\in C^{\infty }\left( \gamma \right) $\ be a primitive of $%
\left( Nu_{\ell }\right) \tau ^{\ast }$, $0\leqslant \ell \leqslant N$.
Assume that $h=\left( u_{\ell }+iv_{\ell }\right) _{0\leqslant \ell
\leqslant N}$\ is an embedding of $\gamma $\ into $\mathbb{C}^{N+1}$. Then (%
\ref{F/ mom affine}) is a necessary condition to the existence of a Riemann
surface $\mathcal{X}$ with almost smooth boundary $\gamma $ such that each $%
u_{\ell }$\ extends to $\mathcal{X}$\ as a harmonic function with harmonic
conjugate function. The converse is true when $\gamma $ is connected and
suitably oriented.\vspace{-3pt}
\end{proposition}

\noindent \textbf{Remarks. }Of course, the above conclusion means that $%
\mathcal{X}$ is a Riemann surface where each $h_{\ell }$ extends
holomorphically.\vspace{-3pt}

\begin{proof}
If $\mathcal{X}$ exists with the required properties, $\delta =h\left(
\gamma \right) $ bounds ,in the sense of current, $h\left( \mathcal{X}%
\right) \backslash \delta $ which is a complex curve of finite volume of $%
\mathbb{C}^{n}\backslash \delta $. Cauchy theorem implies then that (\ref{F/
mom affine}) is verified. If (\ref{F/ mom affine}) is satisfied, \cite%
{HaR-LaB1975} produces a holomorphic 1-chain $Y$ such that $dY=\left[ \delta %
\right] $ ; since $\gamma $ is connected, $Y=\left[ \mathcal{Y}\right] $ for
a suitable orientation of $\gamma $. A normalization of $\overline{\mathcal{Y%
}}$ constructed as in the proof of th.~\ref{T/ caract - a} gives a suitable $%
\mathcal{X}$.
\end{proof}

\end{document}